\documentclass[10pt]{article}
\usepackage{amssymb,amsmath,color,graphicx}
\usepackage{color}
\usepackage{fancybox}
\usepackage[colorlinks]{hyperref}

\newtheorem{definition}{Definition}[section]
\newtheorem{lemma}{Lemma}[section]
\newtheorem{theorem}{Theorem}[section]
\newtheorem{remark}{Remark}[section]

\newtheorem{proposition}{Proposition}[section]
\newtheorem{corollary}{Corollary}[section]

\def\R{\mathbb R}
\def\Z{\mathbb Z}
\def\N{\mathbb N}
\def\epsilon{\varepsilon}
\def\ds{\displaystyle}

\newcommand\mi{\ensuremath{\ds{-\hskip-3.69mm\int}}}

\allowdisplaybreaks

\begin{document}

\title{The speed of propagation for KPP reaction-diffusion equations within  large drift}
\author{Mohammad El Smaily$^\text{ a}$ and St\'ephane Kirsch$^\text{ b}$ \footnote{Both authors are partially supported by a PIMS postdoctoral fellowship.}
 \footnote{ During the preparation of this work, both authors were partially supported by an NSERC grant under the supervision of Professor Nassif Ghoussoub.}\\
{\footnotesize $^\text{a}$ Department of Mathematical Sciences}\\
{\footnotesize{Carnegie Mellon University, Pittsburgh, PA, 15207, USA}}\\
{\footnotesize$^\text{ b}$  Department of Mathematics, University of British Columbia }
 \\{\footnotesize  $\&$ Pacific Institute for the Mathematical Sciences}
 \\{\footnotesize 1984 Mathematics Road, V6T 1Z2, Vancouver, BC, Canada}\\
{\footnotesize {\tt elsmaily@andrew.cmu.edu} , {\tt stephane.kirsch@ac-guyane.fr}}}

\date{}

\maketitle

\begin{abstract} This paper is devoted to the study of the asymptotic behaviors of the minimal speed of propagation of pulsating travelling fronts solving the Fisher-KPP reaction-advection-diffusion equation within either a large drift, a mixture of large drift and small reaction, or a mixture of large drift and large diffusion. We consider a periodic heterogenous framework and we use the formula of Berestycki, Hamel and Nadirashvili \cite{bhn2} for the minimal speed of propagation to prove the asymptotics in any space dimension $N.$ We express the limits as the maxima of certain variational quantities over the family of ``first integrals'' of the advection field. Then, we perform a
detailed study in the case $N=2$ which leads to a necessary and sufficient condition for the positivity of the asymptotic limit of the minimal speed within a large drift.
\end{abstract}


\textbf{AMS Subjet Classifications:} 35B30, 35K55, 35K57,  35Q80, 35Q92, 37C10, 80A32.

\section{Introduction and main results}\label{intro} In this paper, we study the asymptotics of the minimal speed of propagation of pulsating travelling fronts in the presence of a large incompressible advection field. We consider a reaction-advection-diffusion equation
\begin{equation}\label{het eq}
 \left\{
      \begin{array}{l}
 u_t =\nabla\cdot(A(z)\nabla u)\;+M\,q(z)\cdot\nabla u+f(z,u),\; t\in\,\mathbb{R},\;z\in\,\Omega,
\vspace{3pt} \\
        \nu \cdot A\nabla u =0\;\hbox{ on }\mathbb{R}\times\partial\Omega,
      \end{array}
    \right.
\end{equation}
where $\nu$ stands
for the unit outward normal on $\partial\Omega$ whenever it is nonempty.

The domain $\Omega$ is
$C^3$ nonempty connected open subset of $\mathbb{R}^N$ such that for some integer $1\leq d\leq N,$ and for some $L_1,\cdots,L_d$ positive real numbers, we have
\begin{eqnarray}\label{comega}
    \left\{
      \begin{array}{l}
        \exists\,R\geq0\,;\forall\,(x,y)\,\in\,\Omega\subseteq\R^{d}\times\R^{N-d},\,|y|\,\leq\,R, \\
        \forall\,(k_1,\cdots,k_d)\in\,L_1\mathbb{Z}\times\cdots\,\times L_d\mathbb{Z},
        \quad\displaystyle{\Omega\;=\;\Omega+\sum^{d}_{k=1}k_ie_i},
      \end{array}
    \right.
\end{eqnarray}
where $\;(e_i)_{1\leq i\leq N}\;$ is the canonical basis of
$\mathbb{R}^N.$ In other words, $\Omega$ is bounded in the $y-$direction and periodic in $x.$ As archetypes of the domain $\Omega,$ we may have the whole space $\R^{N}$ which corresponds for $d=N$ and $L_1,\cdots, L_N$ any array of positive real numbers. We may also have the whole space $\R^N$ with a periodic
array of holes or an infinite cylinder with an oscillating boundary.
In this periodic situation, we call
\begin{equation}\label{periodicity cell}
C=\{(x,y)\in\,\Omega
;\; x_{1}\in(0,L_1),\cdots,x_{d}\in(0,L_d)\}
\end{equation}  the
periodicity cell of $\Omega.$ We also give the following definition:
\begin{definition}[$L$-periodic fields] A field
$w:\Omega\rightarrow\,\mathbb{R}^N$ is said to be $L$-periodic with
respect to $x$
 if $w(x_1+k_1,\cdots,x_d+k_d\,,y)=w(x_1,\cdots,x_d,y)$ almost everywhere in
$ \Omega,$ and for all
$\displaystyle{k=(k_1,\cdots\,,k_d)\in\prod^{d}_{i=1}L_i\mathbb{Z}}.$
\end{definition}

The diffusion matrix $A(x,y)=(A_{ij}(x,y))_{1\leq i,j\leq
 N}$ is a \textit{symmetric} $C^{2,\delta}(\,\overline{\Omega}\,)$ (with $\delta
\,>\,0$) matrix field satisfying
\begin{eqnarray}\label{cA}
    \left\{
      \begin{array}{l}
        A\; \hbox{is $L$-periodic with respect to}\;x, \vspace{3 pt}\\
        \exists\,0<\alpha_1\leq\alpha_2,\forall(x,y)\;\in\;\Omega,\forall\,\xi\,\in\,\mathbb{R}^N,\vspace{3 pt}\\
       \displaystyle{ \alpha_1|\xi|^2 \;\leq\;\sum_{1\leq i,j\leq N}\,A_{ij}(x,y)\xi_i\xi_j\,\;\leq\alpha_2|\xi|^2.}
      \end{array}
    \right.
\end{eqnarray}

The underlying advection $q(x,y)=(q_1(x,y),\cdots,q_N(x,y))$  is a
$C^{1,\delta}(\overline{\Omega})$ (with $\delta>0$) vector field
satisfying
\begin{equation}\label{cq}
    \left\{
      \begin{array}{ll}
        q\quad\hbox{is $L$- periodic with respect to }\;x, & \hbox{} \\
        \nabla\cdot q=0\quad\hbox{in}\; \overline{\Omega}, \\
         q\cdot\nu=0\quad \hbox{on}\;\partial\Omega\hbox{ (when $\partial\Omega\neq\emptyset$)}, \\
        \forall\,1\leq i\leq d, \quad\displaystyle{\int_{C}q_i\;dx\,dy =0}\hbox{.}
      \end{array}
    \right.
\end{equation}

Concerning the nonlinearity  $f=f(x,y,u),$ it is a nonnegative
function defined in $\overline{\Omega}\,\times[0,1],\;$ such that
\begin{eqnarray}\label{cf1}
    \left\{
      \begin{array}{ll}
        f\geq0, f\;\hbox{ is $L$-periodic with respect to }\; x, \hbox{ and of class }C^{1,\delta}(\overline{\Omega}\times[0,1]),\vspace{3 pt}\\
        \forall\,(x,y)\in\,\overline{\Omega},\quad \displaystyle{f(x,y,0)=f(x,y,1)=0 } \hbox{,} \vspace{3 pt}\\
        \exists \,\rho\in(0,1),\;\forall(x,y)\,\in\overline{\Omega},\;\displaystyle{\forall\, 1-\rho\leq s \leq
s'\leq1,}\;
\displaystyle{f(x,y,s)\;\geq\,f(x,y,s')} \hbox{,} \vspace{3 pt}\\
        \forall\,s\in(0,1),\; \exists \,(x,y)\in\overline{\Omega}\;\hbox{ such that }\;f(x,y,s)>0  \hbox{,} \\
        \forall\,(x,y)\in\overline{\Omega},\quad \zeta(x,y):=\displaystyle{f'_{u}(x,y,0)=\lim_{u\rightarrow\,0^+}\frac{f(x,y,u)}{u}>0}  \hbox{,}
      \end{array}
    \right.
\end{eqnarray}
 with the additional ``KPP'' assumption (referring  to \cite{KPP} by Kolmogorov, Petrovsky and  Piskunov)
 \begin{equation}\label{cf2}
 \forall\, (x,y,s)\in\overline{\Omega}\times(0,1),~0<f(x,y,s)\leq f'_u(x,y,0)\times\,s.
\end{equation}
An archetype of $f$ is $(x,y,u)\mapsto u(1-u)h(x,y)$ defined on $\overline{\Omega}\times[0,1]$ where $h$ is a positive $C^{1,\delta}(\,\overline{\Omega}\,)$ $L$-periodic function.

In all of this paper, $e\in\R^{d}$ is a fixed unit vector and $\tilde{e}:=(e,0,\cdots,0)\in\R^{N}.$ A pulsating travelling front propagating in the direction of $-e$ within a speed $c\neq0$ is a solution $u=u(t,x,y)$ of (\ref{het eq}) for which there exists a function $\phi$ such that
$u(t,x,y)=\phi(x\cdot e+ct,x,y),$  $\phi$ is $L$-periodic in $x$ and $$\lim_{s\rightarrow-\infty}\phi(s,x,y)=0\hbox{ and }\lim_{s\rightarrow+\infty}\phi(s,x,y)=1,$$ uniformly in $(x,y)\in\overline{\Omega}.$

In the same setting as in this paper, it was proved in \cite{bh} and \cite{bhn2} that for all $\Omega,$ $A,\,q,$ and $f$ satisfying (\ref{comega}), (\ref{cA}), (\ref{cq}), (\ref{cf1}-\ref{cf2}) respectively, there exists $\ds{c^{*}_{\Omega,A,q,f}(e)},$ called the \emph{minimal speed of propagation}, such that pulsating travelling fronts exist for $c\geq \ds{c^{*}_{\Omega,A,q,f}(e)}.$ This result extended that of \cite{KPP} which proved that $c^*(e)=2\sqrt{f'(0)}$ in a ``homogenous'' framework where $f=f(u),$ $A=I_N$ (the identity matric) and there is no advection $q.$ A variational formula for the minimal speed $\ds{c^{*}_{\Omega,A,q,f}(e)}$ involving the principal eigenvalue of an elliptic operator was proved in \cite{bhn2} and \cite{weinberger1}. Moreover,  El Smaily \cite{El Smaily min max} proved a $\min$-$\max$ formula for the minimal speed. Many asymptotic behaviors of the minimal speed within large or small diffusion and reaction coefficients and many homogenized speeds were found in \cite{El Smaily} and \cite{EHR}. In \cite{El Smaily}, we have the asymptotic behavior of the minimal speed within a mixture of large diffusion and large advection. Precisely, in Theorem 4.1 of \cite{El Smaily}, it was proved that for all $0\leq\gamma\leq 1/2$ and under the condition $\nabla\cdot A\tilde{e}\equiv 0$ in $\Omega,$
\begin{equation}\label{large diffusion and advection}
\displaystyle{\lim_{M\rightarrow+\infty}\frac{c_{\Omega,MA,\displaystyle{M^{\gamma}\,q},f}^*(e)}{\sqrt{M}}=2\sqrt{\mi_{\!\!\!\!_C}\tilde{e}\cdot A\tilde{e}(x,y)dx\,dy}\,\sqrt{\mi_{\!\!\!\!_C}\zeta(x,y)dx\,dy}},
\end{equation}
where $\zeta(x,y)$ is given in (\ref{cf1}).

 In this paper, we are interested in the asymptotic behavior of ${\ds{c^{*}_{\Omega,A,M\,q,f}(e)}}/{M}$ as $M\rightarrow+\infty$ and also the asymptotic behaviors of the minimal speed within a mixture of large advection and  small reaction or large diffusion. In \cite{bhn1}, it has been proved that
$\liminf_{M\rightarrow+\infty}\ds{\frac{\ds{c^{*}_{\Omega,A,M\,q,f}(e)}}{M}}$ (resp. $\limsup$) are finite. Upper and lower bounds of these $\liminf$
and $\limsup$ were given in \cite{bhn1} in terms of the ``first integrals'' of the advection field $q.$ The family of first integrals of $q$ and two corresponding sub-families will be used in this paper and we recall their definitions below. In Heinze \cite{Heinze convection}, the limit was given in the case of shear flows $q=(q_{_{1}}(y),0,\ldots,0)$ where $e=(1,0,\cdots,0).$ An interesting result about the existence of the limit of $c^{*}_{\Omega,A,M\,q,f}(e)$ as $M\rightarrow+\infty$ (where $\Omega=\R^N$) and several examples of the advection field in 2D and 3D were given in Corollary 1.3 and Section 3 of \cite{Zlatos Ryzhik}.
\begin{definition}[First integrals]\label{first integral} The family of first integrals of $q$ is defined by
\begin{equation*}
\begin{array}{ll}
\mathcal{I}:=&\left\{w\in H^{1}_{loc}(\Omega),\,w\neq0,\;w \hbox{ is } L-\hbox{periodic in $x,$ and }\right. \vspace{4 pt} \\
&\left.q\cdot\nabla w=0\hbox{ almost everywhere in }\Omega\right\}.
\end{array}
\end{equation*}
Having a matrix $A=A(x,y)$ of the type (\ref{cA}), we also define
\begin{equation}\label{I1}
\mathcal{I}_1^A:=\left\{w\in\,\mathcal{I},\hbox{ such that }\int_C\zeta w^2\geq\int_C\nabla w\cdot A\nabla w\right\},
\end{equation}
and
\begin{equation*}
\mathcal{I}_2^A:=\left\{w\in\,\mathcal{I},\hbox{ such that }\int_C\zeta w^2\leq\int_C\nabla w\cdot A\nabla w\right\}.
\end{equation*}
\end{definition}

\begin{remark}[more about $\mathcal{I}$]\label{about I}
The set $\mathcal{I}\cup \{0\}$ is a closed subspace of $H^{1}_{loc}(\Omega).$ Moreover, one can see that if $w\in\mathcal{I}$ is a first integral of $q$
and  $\eta:\R\rightarrow\R$ is a Lipschitz function, then $\eta\circ w\in\mathcal{I}.$
\end{remark}

The following theorem gives the asymptotic behavior of the minimal speed in the presence of a large advection:
\begin{theorem}\label{main}
We fix a unit direction $e\in\R^d$ and assume that the diffusion matrix $A$ and the nonlinearity $f$ satisfy (\ref{cA}), (\ref{cf1}) and (\ref{cf2}).
Let $q$ be an advection field which satisfies (\ref{cq}). Then,
\begin{equation}\label{large advection}
\lim_{M\rightarrow+\infty}\ds{\frac{\ds{c^{*}_{\Omega,A,M\,q,f}(e)}}{M}}=\ds\max_{\ds{w\in \mathcal{I}_1^A}}\frac{\ds\int_C(q\cdot\tilde{e})\,w^2}{\ds\int_C w^2}.
\end{equation}
\end{theorem}
The proof of this theorem will be done later in Section \ref{proof}.
\begin{remark}
It is worth mentioning that the presence of a large advection $M^\gamma q$ ($0\leq\gamma\leq1/2$) has no influence on the limit in (\ref{large diffusion and advection}) whenever a large diffusion $M A$ applies. However, the limit in (\ref{large advection}) depends on $A$ and $f$ via the set $\mathcal{I}_1^A$ and explicitly on the advection $q.$
\end{remark}

\begin{theorem}[large advection with small reaction or large diffusion]\label{mixing large advection with diff and reac}
Assume that $\Omega,A,q$ and $f$ satisfy (\ref{comega}), (\ref{cA}), (\ref{cq}) and (\ref{cf1}-\ref{cf2}) respectively.  Let $e\in\R^d$ be any unit direction. For any $\epsilon>0,$ $B>0$ we call $\ds{c^{*}_{\Omega,A,M\,q,\epsilon f}(e)}$ {\rm(}resp. $\ds{c^{*}_{\Omega,B\,A,M\,q,f}(e)}${\rm)} the minimal speed of propagation (in the direction of $-e$)
of the reaction-advection-diffusion equation with an advection field $Mq$ and a reaction term $\epsilon f$ (resp. diffusion term $B\,A$). Then,
\begin{equation}\label{large advection small reaction}
\lim_{\epsilon\rightarrow0^+}\lim_{M\rightarrow+\infty}\ds{\frac{\ds{c^{*}_{\Omega,A,M\,q,\epsilon f}(e)}}{M\sqrt{\epsilon}}}=\ds\left(\,2\;\frac{\sqrt{\int_{C}\zeta}}{|C|}\,\right)\ds\max_{\ds{w\in \mathcal {I}}}\frac{ \ds\int_C(q\cdot\tilde{e})\,w}{\ds\sqrt{\int_C \nabla w\cdot A\nabla w}},
\end{equation}
and
\begin{equation}\label{large advection large diffusion}
\lim_{B\rightarrow+\infty}\lim_{M\rightarrow+\infty}\ds{\frac{\ds{c^{*}_{\Omega,B\,A,M\,q,f}(e)\times \sqrt{B}}}{M}}=\left(\,2\;\frac{\sqrt{\int_{C}\zeta}}{|C|}\,\right)\ds\max_{\ds{w\in \mathcal {I}}}\frac{ \ds\int_C(q\cdot\tilde{e})\,w}{\ds\sqrt{\int_C \nabla w\cdot A\nabla w}}.
\end{equation}
\end{theorem}
The proof of this theorem will be done in Subsection \ref{proof in mixed cases} below. We mention that many difficulties arose, while demonstrating this result, due to the consideration of a \emph{ heterogeneous framework}. Roughly speaking, the fact that the growth $f'_u(x,y,0)=\zeta(x,y)$ and  the diffusion $A=A(x,y)$ depend on space variables creates a difficulty in choosing a maximizer of the right hand side of (\ref{large advection small reaction}) which should satisfy many properties (see Step 4 of the proof for details). We mention that the above result was proved in the \emph{homogeneous case} ($\zeta=f'(0)$ and $A=Id$) by Zlato\v{s} \cite{zlatos}. In the present paper, we will give the proof of these asymptotics in a general framework.
\vskip0.35cm
Furthermore, in Section \ref{N=2 section} of this paper, we will give more details about the family of first integrals $\mathcal{I}$ and about integrals of the form   $\ds\int_C(q\cdot\tilde{e})\,w^2$ (where $w\in\mathcal{I}$) in the case where $N=2$. This will give  necessary and sufficient conditions, expressed in terms of the nature of the drift $q,$ for the limit (\ref{large advection}) to be null or not. In this context, we have Theorem \ref{case N=2} which will be announced after the following definition.

\begin{definition}\label{trajectory def}
Assume that $N=2$ and that $\Omega$ and $q$ satisfy (\ref{comega}) and (\ref{cq}). Let $x \in \Omega$ such that $q(x) \neq 0.$ The trajectory of $q$ at $x$ is the largest (in the sense of inclusion) connected differentiable curve $T(x)$ in $\Omega$ verifying:

(i) $x \in T(x)$,

(ii) $\forall y \in T(x)$, $q(y) \neq 0$,

(iii)  $\forall y \in T(x)$, $q(y)$ is tangent to $T(x)$ at the point $y$.
\end{definition}
In the following lemma, we describe the family of ``unbounded periodic trajectories'' of a vector field $q$. The proof of this lemma will be done in Section \ref{N=2 section}.

\begin{lemma}[unbounded periodic trajectories]\label{periodic trajs}
Let $T(x)$ be an unbounded periodic trajectory of $q$ in $\Omega,$ that is there exists $\mathbf{a} \in L_1\Z\times L_2\Z\setminus\{0\}$ (resp. $L_1\Z\times \{0\}\setminus\{0\}$) when $d=2$ (resp. $d=1$) such that $T(x)=T(x)+\mathbf{a}.$ In this case, we say that $T(x)$ is $\mathbf{a}-$periodic. Then, if $T(y)$ is another unbounded periodic trajectory of $q,$ $T(y)$ is also $\mathbf{a}-$periodic.

Moreover, in the case $d=1,$ $\mathbf{a}=L_1e_1.$ That is, all the unbounded periodic trajectories of $q$ in $\Omega$ are $L_1e_1-$periodic.
\end{lemma}

\begin{theorem}\label{case N=2}
Assume that $N=2$ and that $\Omega$ and $q$ satisfy (\ref{comega}) and (\ref{cq}) respectively. The two following statements are equivalent:

(i) There exists $w \in \mathcal{I},$ such that $\displaystyle \int_C qw^2 \neq 0$.

(ii) There exists a periodic unbounded trajectory $T(x)$ of $q$ in $\Omega$.

Moreover, if (ii) is verified and $T(x)$ is $\mathbf{a}-$periodic, then for any $w\in\mathcal{I}$ we have $\ds{\int_Cq\,w^2\in\R\mathbf{a}}.$
\end{theorem}

\begin{remark}
The periodicity assumption on the trajectory in (ii) is crucial. Indeed there may exist unbounded trajectories which are not periodic, even though the vector field $q$ is periodic. Consider the following function $\phi$:
$$
\phi(x,y) :=\left\{\begin{array}{l}\ds{ e^{-\frac{1}{\sin^2(\pi y)}}\sin(2\pi(x+\ln(y-[y])))\text{ if }y\not\in\mathbb{Z}},\vspace{3pt}\\
0\text{ otherwise},
\end{array}\right.
$$
where $[y]$ denotes the integer part of $y$. This function is $C^\infty$ on $\mathbb{R}^2$, and $1$-periodic in $x$ and $y$. Hence the vector field
$$
q = \nabla^\perp \phi
$$
is also $C^\infty$, $1$-periodic in $x$ and $y$, and verifies $\int_{[0,1]\times[0,1]} q = 0$ with $\nabla\cdot q\equiv0.$ A quick study of this vector field shows that the part of the graph of $x \mapsto e^{-x}$ lying between $y=0$ and $y=1$ is a trajectory of $q$, and is obviously unbounded and not periodic. However there exist no periodic unbounded trajectory for this vector field, so the theorem asserts that  for all $w \in \mathcal{I}$ we have
$$
\int_C qw^2 =0.
$$
\end{remark}

As a direct consequence of Theorem \ref{main} and Theorem \ref{case N=2}, we get the following corollary about the asymptotic behavior of the minimal speed within large drift:

\begin{corollary}\label{NSC for lim to be zero} Assume that $N=2$ and that $\Omega,$ $A,$ $q$ and $f$ satisfy the conditions (\ref{comega}), (\ref{cA}), (\ref{cq}) and (\ref{cf1}-\ref{cf2}) respectively. Then,

(i) If there exists no periodic unbounded trajectory of $q$ in $\Omega,$ then
$$
\lim_{M\rightarrow+\infty}\ds{\frac{\ds{c^{*}_{\Omega,A,M\,q,f}(e)}}{M}}=0,
$$ for any unit direction $e.$

(ii) If there exists a periodic unbounded trajectory $T(x)$ of $q$ in $\Omega$ {\rm(}which will be $\mathbf{a}-$periodic for some vector $\mathbf{a}\in\mathbb{R}^2${\rm)} then
\begin{equation}\label{null equivalence}
\lim_{M\rightarrow+\infty}\ds{\frac{\ds{c^{*}_{\Omega,A,M\,q,f}(e)}}{M}} >0\,\Longleftrightarrow\,\tilde{e}\cdot\mathbf{a}\neq0.
\end{equation}
We mention that in the case where $d=1,$ we have $\tilde{e}=\pm e_1.$  Lemma \ref{periodic trajs} yields that $\tilde{e}\cdot\mathbf{a}=\pm L_1\neq0.$ Referring to (\ref{null equivalence}), we can then write,  for $d=1,$
\begin{equation}\label{null equivalence d=1}
\lim_{M\rightarrow+\infty}\ds{\frac{\ds{c^{*}_{\Omega,A,M\,q,f}(e)}}{M}} >0\,\Longleftrightarrow\,\text{(there exists a periodic unbounded trajectory $T(x)$ of $q$ in $\Omega$)}.
\end{equation}

\end{corollary}
It is worth mentioning that, in the above corollary, the conditions for which the limit is null or not are expressed only in terms the advection field $q$ and, moreover,
it is easy to check if they are verified by $q$ or not.
\begin{remark}
In (ii), the simplest example is when $q$ is a shear flow (i.e. $q(x_1,x_2)=(q_1(x_2),0)$). In that case, the limit (\ref{large advection}) is positive if and only if $\tilde{e}$ is not perpendicular to the flow lines of $q$ (this condition means that the first component of $\tilde{e}$ is not zero).
\end{remark}

\subsection{Outline of the rest of the paper} After the announcement of the main results in Section \ref{intro}, we are going to
prove, in Section \ref{proof}, the asymptotics of the minimal speed within large drift in any space dimension $N$. This section will be divided into
two subsections. In the first one, we prove  (\ref{large advection}) which deals with the asymptotic behavior of the speed in the presence of large advection \textbf{only}  and then, in Subsection \ref{proof in mixed cases}, we prove  (\ref{large advection small reaction}) which concerns the asymptotic behavior of the speed in a mixture of large drift and small reaction (or a mixture of large drift and a large diffusion). In  Section \ref{N=2 section}, we prove many auxiliary lemmas which  lead to  the proof of Theorem \ref{case N=2} and Corollary \ref{NSC for lim to be zero} in the case where $N=2.$

\section{Proofs of the asymptotic behaviors in any dimension $N$}\label{proof}
  Theorems \ref{main} and \ref{mixing large advection with diff and reac} were announced for domains $\Omega\subseteq\R^N $
where $N$ could be any dimension. We divide the present section into two subsections. The first subsection deals with the case where we have only a large advection and the second one
deals with the case where we have large advection mixed with a small reaction or a large diffusion.

\subsection{Case of large advection (proof of Theorem \ref{main})}\label{proof in case large drift}
Here we prove Theorem \ref{main}. For this, we start by a proposition which will play an important role
in the proof. \textbf{For the sake of simplicity, we suppose that the diffusion matrix}  $A=A(x,y)=Id$ where $Id$ is the identity matrix of $M_N(\R).$ In the case of any matrix $A$ satisfying (\ref{cA}), the proof of (\ref{large advection})  is very similar to that in the case $A=Id.$ Indeed, we point out the simple differences in Remark \ref{for general diffusion A} below.

 Over all this proof,
 since the only parameter in $\ds{c^{*}_{\Omega,A,M\,q,f}(e)}$  is the factor $M$ in front of the  advection $q,$ we will write
$$\ds{c^{*}(M):=c^{*}_{\Omega,A,M\,q,f}(e)}.$$
After supposing that $A=Id,$ the subsets $\mathcal{I}_1^A$ and $\mathcal{I}_2^A$ are respectively given by
\begin{equation}\label{I1 A=id}
\mathcal{I}_1:=\{w\in\,\mathcal{I},\hbox{ such that }\int_C\zeta w^2\geq\int_C|\nabla w|^2\},
\end{equation}
and
\begin{equation}\label{I2 A=id}
\mathcal{I}_2:=\{w\in\,\mathcal{I},\hbox{ such that }\int_C\zeta w^2\leq\int_C|\nabla w|^2\}.
\end{equation}

\begin{definition}
We define $g: [0,+\infty) \rightarrow \R$ by
\begin{equation} \label{defg}
g(\lambda) := \sup_{w \in \mathcal{I}} \frac{1}{\int_C w ^2}\left[\int_C \left(\zeta w^2 - |\nabla w|^2 \right) + \lambda \int_C (q \cdot \tilde{e})w^2\right],
\end{equation}
and we define $h: (0,+\infty) \rightarrow \R$ by
\begin{equation} \label{defh}
h(\lambda) := \frac{g(\lambda)}{\lambda}.
\end{equation}
\end{definition}

\begin{remark}\label{remsup}
We can replace the supremum by a maximum in (\ref{defg}). Indeed, consider for a fixed $\lambda>0$ a maximizing sequence $\{w_n\}_n$ with $\|w_n\|_{L^2(C)}=1$. We have
$$
\left[\int_C \left(\zeta w_n^2 - |\nabla w_n|^2 \right) + \lambda \int_C (q \cdot \tilde{e})w_n^2\right] \xrightarrow[n \rightarrow +\infty]{} g(\lambda).
$$
The sequence $\{w_n\}_n$ is then bounded in $H^1_{loc}(\Omega)$ and we can extract a subsequence converging weakly in $H^1_{loc}(\Omega)$ and strongly in $L^2_{loc}(\Omega)$ to $w_0$. We then have  $w_0 \in \mathcal{I}$, $\|w_0\|_{L^2(C)}=1.$ Since $\displaystyle \liminf_{n \rightarrow \infty} \int_C |\nabla w_n|^2 \geq \int_C |\nabla w_0|^2,$ we then get
$$
\left[\int_C \left(\zeta w_0^2 - |\nabla w_0|^2 \right) + \lambda \int_C (q \cdot \tilde{e})w_0^2\right] \geq g(\lambda),
$$
and by definition of the supremum, the previous inequality is an equality, which means by the way that the weak convergence in $H^1_{loc}(\Omega)$ is in fact a strong convergence.
\end{remark}

\begin{proposition}\label{prop}
The functions $g$ and $h$ satisfy the following properties:

(i) The function $g$ is convex on $[0,+\infty),$ and moreover,  $g$ and $h$ are continuous on their domains and take values in $(0,+\infty)$.\\

(ii) $ \displaystyle h(\lambda) \xrightarrow[\lambda \rightarrow +\infty]{} \sup_{w \in \mathcal{I}} \frac{\int_C (q \cdot \tilde{e})w^2 }{\int_C w^2}$.\\

(iii) Either $h$ is convex and decreasing on $(0,+\infty)$ or $h$ attains a global minimum at some point $\lambda_0 >0$.\\

(iv) If $h$ is convex decreasing on $(0,+\infty)$, then we have

$$
h(\lambda) \xrightarrow[\lambda \rightarrow +\infty]{} \max_{w \in \mathcal{I}_1} \frac{\int_C (q \cdot \tilde{e})w^2 }{\int_C w^2} = \sup_{w \in \mathcal{I}} \frac{\int_C (q \cdot \tilde{e})w^2 }{\int_C w^2}.
$$\\

(v) If $h$ attains its minimum at $\lambda_0>0$, then we have
\begin{equation} \label{min}
h(\lambda_0) = \max_{w \in \mathcal{I}_1} \frac{\int_C (q \cdot \tilde{e})w^2 }{\int_C w^2}.
\end{equation}
\end{proposition}

{\bf Proof of (i):} $g$ is the supremum of affine functions, so it is convex, and hence continuous. Since $\int_C \zeta >0$, and since the constant functions belong to $\mathcal{I},$ we have
$$
\forall \lambda \geq 0,~~g(\lambda) \geq \frac{\int_C \zeta}{|C|} >0.
$$
Hence,  $g(\lambda) >0$ for any $\lambda \geq 0.$ Besides,  $h$ is well defined and continuous, and $h(\lambda) >0$ for any $\lambda >0$.
\vskip 0.3cm

{\bf Proof of (ii):} For each $k \in \mathbb{N},$ we define
$$
\mathcal{I}^k:=\left\{w \in \mathcal{I}, \text{ such that } \int_C(\zeta w^2 - |\nabla w|^2) \geq -k \int_C w^2 \right\}
$$
and
$$
\ds{h_k(\lambda) := \sup_{w \in \mathcal{I}^k} \frac{1}{\int_C w ^2}\left[\frac{1}{\lambda}\int_C \left(\zeta w^2 - |\nabla w|^2 \right) + \int_C (q \cdot \tilde{e})w^2\right]}.
$$
Obviously, $h_k(\lambda) \leq h(\lambda)$ for any $\lambda >0$ because the supremum is taken over a smaller set. Moreover, a simple computation gives
$$
h_k(\lambda) \xrightarrow[\lambda \to +\infty]{} \sup_{w \in \mathcal{I}^k} \frac{\int_C (q \cdot \tilde{e})w^2 }{\int_C w^2}.
$$
Hence, for every $k \in \mathbb{N},$ we have
$$
\liminf_{\lambda \to +\infty} h(\lambda) \geq \sup_{w \in \mathcal{I}^k} \frac{\int_C (q \cdot \tilde{e})w^2 }{\int_C w^2},
$$
and since $\displaystyle \mathcal{I} = \bigcup_{k \in \mathbb{N}} \mathcal{I}^k,$ we get
$$
\liminf_{\lambda \to +\infty} h(\lambda) \geq \sup_{w \in \mathcal{I}} \frac{\int_C (q \cdot \tilde{e})w^2 }{\int_C w^2}.
$$
On the other hand,
$$
h(\lambda) \leq \frac{\|\zeta\|_\infty}{\lambda} + \sup_{w \in \mathcal{I}} \frac{\int_C (q \cdot \tilde{e})w^2 }{\int_C w^2},
$$
which gives
$$
\limsup_{\lambda \to +\infty} h(\lambda) \leq \sup_{w \in \mathcal{I}} \frac{\int_C (q \cdot \tilde{e})w^2 }{\int_C w^2},
$$
and completes the proof of (ii).
\vskip 0.3cm

{\bf Proof of (iii):} We know from (ii) that $h(\lambda)$ converges when $\lambda \rightarrow +\infty$. Moreover, since $g(0)>0$ we then have $h(\lambda) \rightarrow +\infty$ as $\lambda \rightarrow 0.$  We distinguish now two different cases:

{\bf Case 1:} Suppose that for any $ \lambda >0,$ we have $\displaystyle h(\lambda) > \lim_{\lambda\rightarrow+\infty}h(\lambda)$. Thus, for a fixed $\lambda>0,$ the definition of the limit yields the existence of $\lambda_1>\lambda$ such that $h(\lambda)>h(\lambda_1)$. Let then $w$ be such that $\|w\|_{L^2(C)}$=1 and
$$
h(\lambda) = \frac{\int_C(\zeta w^2 - |\nabla w|^2)}{\lambda} + \int_C(q \cdot \tilde{e})w^2.
$$ (the existence of $w$ follows from Remark \ref{remsup}).
 From the definition of $h,$ we can conclude that
$$
\frac{\int_C(\zeta w^2 - |\nabla w|^2)}{\lambda} + \int_C(q \cdot \tilde{e})w^2 = h(\lambda)>h(\lambda_1) \geq \frac{\int_C(\zeta w^2 - |\nabla w|^2)}{\lambda_1} + \int_C(q \cdot \tilde{e})w^2,
$$
which gives
$$
\left(\frac{1}{\lambda} - \frac{1}{\lambda_1} \right) \int_C(\zeta w^2 - |\nabla w|^2) \geq 0.
$$
Having $\lambda < \lambda_1,$ we get
$$
\int_C(\zeta w^2 - |\nabla w|^2) \geq 0.
$$

Thus, for every $\lambda>0$, the maximum in the definition of $h(\lambda)$ is attained in $\mathcal{I}_1$. Therefore,  $h$ can be rewritten in this case as follows
$$
h(\lambda) = \max_{w \in \mathcal{I}_1}\frac{\int_C(\zeta w^2 - |\nabla w|^2)}{\lambda \int_C w^2} + \frac{\int_C(q \cdot \tilde{e})w^2}{\int_C w^2}.
$$
In this formulation of $h,$ we maximize over $\mathcal{I}_1$.  The map $$\lambda\mapsto\frac{\int_C(\zeta w^2 - |\nabla w|^2)}{\lambda \int_C w^2}$$ is convex when $w\in\mathcal{I}_1.$  Hence,
 $h$ is the supremum of convex functions and is then convex. Moreover, $h$ converges when $\lambda \rightarrow +\infty$, and $h(\lambda)> \lim_{+\infty}h$, which, with the convexity of $h$, implies that $h$ is decreasing on $(0,+\infty)$.\\

{\bf Case 2:} There exists $\lambda>0$ such that $\displaystyle h(\lambda) \leq \lim_{\lambda\rightarrow+\infty}h(\lambda)$. By continuity, there exists $\lambda_0>0$ such that
$$
h(\lambda_0) = \min_{\lambda>0} h(\lambda).
$$
\vskip 0.3cm

{\bf Proof of (iv):} In the case where $h$ is convex and decreasing over $(0,+\infty),$ we know that
$$
h(\lambda) = \max_{w \in \mathcal{I}_1}\frac{\int_C(\zeta w^2 - |\nabla w|^2)}{\lambda \int_C w^2} + \frac{\int_C(q \cdot \tilde{e})w^2}{\int_C w^2}.
$$
Thus,
$$
\max_{w \in \mathcal{I}_1} \frac{\int_C(q \cdot \tilde{e})w^2}{\int_C w^2} \leq h(\lambda) \leq \frac{\|\zeta\|_\infty}{\lambda} + \max_{w \in \mathcal{I}_1} \frac{\int_C(q \cdot \tilde{e})w^2}{\int_C w^2}.
$$
We conclude that
$$
h(\lambda) \xrightarrow[\lambda \rightarrow +\infty]{} \max_{w \in \mathcal{I}_1} \frac{\int_C(q \cdot \tilde{e})w^2}{\int_C w^2}.
$$\vskip 0.3cm

{\bf Proof of (v):} We use several claims to prove this last part of Proposition \ref{prop}. The proofs of these claims are postponed at the end.\\ \vskip0.2cm
\noindent{\bf $\bullet$ Claim 1:} There exist $w_1 \in \mathcal{I}_1$ and $w_2 \in \mathcal{I}_2$ such that
$$
h(\lambda_0) = \frac{\int_C(\zeta w_1^2 - |\nabla w_1|^2)}{\lambda_0 \int_C w_1^2} + \frac{\int_C(q \cdot \tilde{e})w_1^2}{\int_C w_1^2}= \frac{\int_C(\zeta w_2^2 - |\nabla w_2|^2)}{\lambda_0 \int_C w_2^2} + \frac{\int_C(q \cdot \tilde{e})w_2^2}{\int_C w_2^2}.
$$
{\bf $\bullet$ Claim 2:} If $w_1 \in \mathcal{I}$ and $w_2 \in \mathcal{I}$ are not proportional and
$$
h(\lambda_0) = \frac{\int_C(\zeta w_1^2 - |\nabla w_1|^2)}{\lambda_0 \int_C w_1^2} + \frac{\int_C(q \cdot \tilde{e})w_1^2}{\int_C w_1^2}= \frac{\int_C(\zeta w_2^2 - |\nabla w_2|^2)}{\lambda_0 \int_C w_2^2} + \frac{\int_C(q \cdot \tilde{e})w_2^2}{\int_C w_2^2},
$$
then for any $0 \leq \theta \leq 1$ and $w_\theta := \theta w_1 + (1-\theta)w_2,$ we have
$$
h(\lambda_0) = \frac{\int_C(\zeta w_\theta^2 - |\nabla w_\theta|^2)}{\lambda_0 \int_C w_\theta^2} + \frac{\int_C(q \cdot \tilde{e})w_\theta^2}{\int_C w_\theta ^2}.
$$
Claim 1 gives us $w_1 \in \mathcal{I}_1$ and $w_2 \in \mathcal{I}_2$ realizing the maximum in the definition of $h(\lambda_0)$. If $w_1$ and $w_2$ are proportional, then $w_1$ (resp. $w_2$) $\in \mathcal{I}_1 \cap \mathcal{I}_2$ and we define $w_0:=w_1$. If not, using claim 2, we know that any convex combination also realizes the maximum in the definition of $h$. By continuity, there exists $\theta_0 \in [0,1]$ such that $w_0 :=w_{\theta_0} \in \mathcal{I}_1 \cap \mathcal{I}_2$ and
$$
h(\lambda_0) = \frac{\int_C(\zeta w_0^2 - |\nabla w_0|^2)}{\lambda_0 \int_C w_0^2} + \frac{\int_C(q \cdot \tilde{e})w_0^2}{\int_C w_0^2},
$$
and since $w_0 \in \mathcal{I}_1 \cap \mathcal{I}_2$, we then have $\int_C(\zeta w_0^2 - |\nabla w_0|^2) = 0$ and consequently,
$$
h(\lambda_0) = \frac{\int_C(q \cdot \tilde{e})w_0^2}{\int_C w_0^2}.
$$
Since $w_0 \in \mathcal{I}_1$, we have
$$
h(\lambda_0) \leq \max_{w \in \mathcal{I}_1} \frac{\int_C(q \cdot \tilde{e})w^2}{\int_C w^2}.
$$
On the other hand, by the definition of $h$ we have
\begin{eqnarray*}
h(\lambda_0) & = & \sup_{w \in \mathcal{I}}\frac{\int_C(\zeta w^2 - |\nabla w|^2)}{\lambda_0 \int_C w^2} + \frac{\int_C(q \cdot \tilde{e})w^2}{\int_C w^2} \\
& \geq & \max_{w \in \mathcal{I}_1} \frac{\int_C(\zeta w^2 - |\nabla w|^2)}{\lambda_0 \int_C w^2} + \frac{\int_C(q \cdot \tilde{e})w^2}{\int_C w^2} \\
& \geq & \max_{w \in \mathcal{I}_1} \frac{\int_C(q \cdot \tilde{e})w^2}{\int_C w^2}.
\end{eqnarray*}
This ends the proof of (v).
\vskip0.35cm
\noindent We are left to prove claims 1 and 2.
\vskip0.5cm

{\bf \underline{Proof of Claim 1}:} Let $\{\lambda_1^p\}_p$ be a sequence such that $\lambda_1^p < \lambda_0$ and $\lambda_1^p \rightarrow \lambda_0$ as $p\rightarrow+\infty.$\\
For each $p\in\N,$ let $w_1^p\in\mathcal{I}$ such that $\int_C (w_1^p)^2 = 1$ and
$$
h(\lambda_1^p) = \frac{1}{\lambda_1^p}\int_C\left(\zeta (w_1^p)^2 - |\nabla w_1^p|^2  \right) + \int_C (q \cdot \tilde{e})(w_1^p)^2.
$$
From the definition of $h$ and owing to the fact that $h(\lambda_1^p) \geq h(\lambda_0),$ we have
\begin{eqnarray*}
h(\lambda_1^p) & = & \frac{1}{\lambda_1^p}\int_C\left(\zeta (w_1^p)^2 - |\nabla w_1^p|^2  \right) + \int_C (q \cdot \tilde{e})(w_1^p)^2 \\
& \geq & h(\lambda_0)  \\
& \geq & \frac{1}{\lambda_0}\int_C\left(\zeta (w_1^p)^2 - |\nabla w_1^p|^2  \right) + \int_C (q \cdot \tilde{e})(w_1^p)^2.
\end{eqnarray*}
However, $\lambda_1^p \leq \lambda_0.$ Thus,
\begin{equation}\label{w_1^P belongs to I1}
\int_C(\zeta (w_1^p)^2 - |\nabla w_1^p|^2) \geq 0,
\end{equation}
which means $w_1^p \in \mathcal{I}_1$. Moreover, (\ref{w_1^P belongs to I1}) yields that $\{ w_1^p \}_p$ is a bounded sequence in $H^1_{loc}(\Omega).$ Therefore, we can extract a subsequence converging weakly in $H^1_{loc}(\Omega)$ and strongly in $L^2_{loc}(\Omega)$ to $w_1 \in \mathcal{I}.$ Since the convergence is strong in $L^2_{loc}(\Omega),$ we get $\int_C (w_1)^2 = 1$. Thanks to the continuity of $h$ with respect to $\lambda,$ we get
\begin{equation}\label{eq1}
\frac{1}{\lambda_1^p}\int_C\left(\zeta (w_1^p)^2 - |\nabla w_1^p|^2  \right) + \int_C (q \cdot \tilde{e})(w_1^p)^2 \xrightarrow[p \rightarrow \infty]{} h(\lambda_0).
\end{equation}
Moreover, we have
\begin{equation} \label{eq2}
h(\lambda_0) \geq \frac{1}{\lambda_0}\int_C\left(\zeta (w_1)^2 - |\nabla w_1|^2  \right) + \int_C (q \cdot \tilde{e})(w_1)^2
\end{equation}
and
\begin{equation} \label{eq3}
\frac{1}{\lambda_1^p}\int_C\left(\zeta (w_1^p)^2\right)+ \int_C (q \cdot \tilde{e})(w_1^p)^2 \xrightarrow[p \rightarrow \infty]{} \frac{1}{\lambda_0}\int_C\left(\zeta (w_1)^2\right)+ \int_C (q \cdot \tilde{e})(w_1)^2.
\end{equation}
The combination of (\ref{eq1}), (\ref{eq2}) and (\ref{eq3}) gives
$$
\limsup_{p \rightarrow \infty}\int_C|\nabla w_1^p|^2 \leq \int_C|\nabla w_1|^2.
$$
On the other hand, the weak convergence $w_1^p \rightharpoonup w_1$ in $H^1_{loc}(\Omega)$ implies that
$$
\liminf_{p \rightarrow \infty}\int_C|\nabla w_1^p|^2 \geq \int_C|\nabla w_1|^2.
$$
Hence, $\{w_1^p\}_p$ converges strongly in $H^1_{loc}(\Omega)$ to $w_1$. We then conclude that $w_1 \in \mathcal{I}_1$ (because $\mathcal{I}_1$ is a closed subset of $H^1_{loc}(\Omega)$) and that
$$
h(\lambda_0) = \frac{\int_C(\zeta w_1^2 - |\nabla w_1|^2)}{\lambda_0 \int_C w_1^2} + \frac{\int_C(q \cdot \tilde{e})w_1^2}{\int_C w_1^2}.
$$
We can use a similar argument (we take $\lambda_2^p > \lambda_0$ such that $\lambda_2^p \rightarrow \lambda_0$ as $p\rightarrow+\infty$ and, for each $p,$ we take $w_2^p$ as a maximizer of $h(\lambda_2^p )$) to get $w_2 \in \mathcal{I}_2$ verifying
$$
h(\lambda_0) = \frac{\int_C(\zeta w_2^2 - |\nabla w_2|^2)}{\lambda_0 \int_C w_2^2} + \frac{\int_C(q \cdot \tilde{e})w_2^2}{\int_C w_2^2}.
$$\\
{\bf \underline{Proof of Claim 2}:} Without loss of generality, we suppose that $\int_C w_1^2 = \int_C w_2^2=1$.
We consider the following functional defined by
$$\forall\,w\in\mathcal{I},~~
E_\lambda(w) :=\frac{\int_C(\zeta w^2-|\nabla w|^2)}{\lambda \int_C w^2} + \frac{\int_C(q \cdot \tilde{e})w^2}{\int_C w^2}.
$$
We have
$$
h(\lambda_0) = \max_{w \in \mathcal{I}} E_{\lambda_0}(w) = E_{\lambda_0}(w_1)
$$
and thereby, $\forall w \in \mathcal{I}$ we have
\begin{equation}\label{deriv}
E'_{\lambda_0}(w_1)w = 0 =  \frac{1}{\lambda_0}\int_C (\zeta w_1 w - \nabla w_1 \cdot \nabla w) + \int_C(q \cdot \tilde{e})w_1\cdot w - h(\lambda_0) \int_C w_1w.
\end{equation}
Now, we compute  $E_{\lambda_0}(w_\theta)$ explicitly.
We have
$$
\left\{
\begin{array}{lll}
\displaystyle \int_C w_\theta^2 & = &\ds{ \theta^2 + (1-\theta)^2 + 2\theta(1-\theta)\int_Cw_1 w_2,} \vspace{3 pt}\\
\displaystyle \int_C \zeta w_\theta^2 & = & \ds{\theta^2 \int_C \zeta w_1^2 + (1-\theta)^2\int_C \zeta w_2^2 + 2\theta(1-\theta)\int_C \zeta w_1 w_2,} \vspace{3 pt}\\
\displaystyle \int_C |\nabla w_\theta|^2 & = & \ds{\theta^2 \ds\int_C |\nabla w_1^2| + (1-\theta)^2 \int_C |\nabla w_2|^2 + 2\theta(1-\theta) \int_C \nabla w_1\cdot \nabla w_2,}\vspace{3 pt}\\
\displaystyle \int_C (q \cdot \tilde{e})w_\theta^2 & = & \ds{\theta^2 \int_C (q \cdot \tilde{e})w_1^2 + (1-\theta)^2\int_C (q \cdot \tilde{e})w_2^2 + 2\theta(1-\theta) \int_C (q \cdot \tilde{e})w_1 w_2},
\end{array}
\right.
$$
and using (\ref{deriv}) with $w = w_2$ we get
$$
E_{\lambda_0}(w_\theta) = \frac{\theta^2h(\lambda_0) +(1-\theta)^2h(\lambda_0) + 2\theta(1-\theta)h(\lambda_0)\int_Cw_1w_2}{\theta^2+(1-\theta)^2 + 2\theta(1-\theta)\int_Cw_1 w_2}.
$$
The denominator is positive because, by assumption, $w_1$ and $w_2$ are not proportional, so we can not have equality in the Cauchy-Schwarz inequality. After simplification, we obtain
$$
E_{\lambda_0}(w_\theta) = h(\lambda_0).
$$
This completes the proof of Proposition \ref{prop}.\hfill$\Box$
\vskip0.45cm
\textbf{Proof of Theorem \ref{main}.}
From the results of \cite{bhn2}, it follows  that for each $M>0,$ the minimal speed $c^*(M)$ is given by
\begin{equation}\label{var formula}
c^*(M)=\min_{\lambda>0}\frac{k(\lambda,M)}{\lambda},
\end{equation}
where $k(\lambda,M)$ is the principal eigenvalue of the elliptic operator $L_\lambda$ defined by
\begin{equation}
L_\lambda \psi:=\displaystyle{\Delta\psi+2\lambda\tilde{e}\cdot
\nabla\psi+M\,q\cdot\nabla\psi}+\displaystyle{[\lambda^2+\lambda M\,
q\cdot\tilde{e}+\zeta]\psi}\hbox{ in }\Omega,
\end{equation}
acting on the set $E_\lambda$
$$
E_\lambda=\left\{\psi=\psi(x,y)\in C^2(\overline{\Omega}), \psi\hbox{ is
$L$-periodic in $x$ and }
\nu\cdot \nabla\psi=-\lambda(\nu\cdot
\tilde{e})\psi\;\hbox{on}\;\partial{\Omega}\right\}.
$$
The principal eigenfunction $\psi^{\lambda,M}$ associated to $k(\lambda,M)$ is positive in $\overline{\Omega}$ and
it is unique up to multiplication by a nonzero real number. The existence of $k(\lambda,M)$ and $\psi^{\lambda,M}$ for any $(\lambda,M)\in\R\times\R,$ and the properties of $k(\lambda,M)$ as a function of $M$ have been studied in \cite{bh} and \cite{bhn2}.
 In particular, the function $\lambda\mapsto k(\lambda,M)$ is convex and $k(\lambda,M)>0$ for all $(\lambda,M)\in(0,+\infty)\times(0,+\infty).$

We want to study the asymptotic behavior of $M\mapsto c^{*}(M)/M$ when $M\rightarrow+\infty.$ For this, we call $\lambda'=\lambda\times M,$
 $\mu(\lambda',M)=k(\lambda,M)$ and $\psi^{\lambda',M}=\psi^{\lambda,M}$ for each
$(\lambda,M)\in(0,+\infty)\times(0,+\infty).$ Referring to formula (\ref{var formula}), we then get
\begin{equation}\label{var form in mu}
\forall \,M>0,~~\frac{c^*(M)}{M}=\min_{\lambda'>0}\frac{\mu(\lambda',M)}{\lambda'}.
\end{equation}
From the properties of $k(\lambda,M),$ we have $\lambda'\mapsto\mu(\lambda',M)$ is convex over $(0,+\infty)$ and $\mu(\lambda',M)>0$ for all
$(\lambda',M)\in (0,+\infty)\times(0,+\infty).$ Moreover, it follows from above that the function $\psi^{\lambda',M}$ and $\mu(\lambda',M)$ are respectively the principal eigenfunction and the principal eigenvalue of the problem
\begin{equation}\label{Leq mu}
\left\{
  \begin{array}{rl}
 \mu(\lambda',M)\psi^{\lambda',M}   =&\displaystyle{\Delta\psi^{\lambda',M}+2\frac{\lambda'}{M}\tilde{e}\cdot
\nabla\psi^{\lambda',M}+M\,q\cdot\nabla\psi^{\lambda',M}}\vspace{3pt}\\
&+\displaystyle{\left[\left(\frac{\lambda'}{M}\right)^2+\lambda'\,
q\cdot\tilde{e}+\zeta\right]\psi^{\lambda',M}}\hbox{ in }\Omega,\vspace{3 pt}\\
\nu\cdot \nabla\psi^{\lambda',M}=&\ds{-\frac{\lambda'}{M}(\nu \cdot
\tilde{e})\psi^{\lambda',M}}\;\hbox{on}\;\partial{\Omega}\hbox{ (whenever }\partial\Omega\neq\emptyset\hbox{)}.
  \end{array}
\right.
\end{equation}

We take any first integral $w\in \mathcal{I}$ of $q,$ we multiply (\ref{Leq mu}) by $\ds{\frac{w^2}{\psi^{\lambda',M}}}$ and integrate by parts
over the periodicity cell $C.$ Using (\ref{cq}), the boundary condition on $\psi^{\lambda',M}$ and the fact that $q\cdot\nabla w=0$ a.e in $\Omega,$  we obtain that
$$\begin{array}{ll}
\ds{\mu(\lambda',M)\int_{C}w^2}=&\ds{\int_{C}\left(\frac{\nabla\psi^{\lambda',M}}{\psi^{\lambda',M}}\,w\right)^2-\frac{\lambda'}{M}\int_{\partial C}\nu\cdot\tilde{e}w^2-2\int_{C}\left(\frac{\nabla\psi^{\lambda',M}}{\psi^{\lambda',M}}\,w\right)\cdot\left(\nabla w-\frac{\lambda'}{M}\tilde{e}w\right)}\vspace{4 pt}\\
&\ds{+\left(\frac{\lambda'}{M}\right)^2\int_C w^2 +\lambda'\int_C q\cdot\tilde{e}\,w^2+\int_{C}\zeta w^2}.
\end{array}
$$
Notice that the boundary term $\frac{\lambda'}{M}\int_{\partial C}\nu\cdot\tilde{e}w^2$ is equal to $\frac{\lambda'}{M}\int_{\partial C}2w\nabla w\cdot \tilde {e}.$
After dividing the previous equation by $\lambda',$ we get
\begin{equation}\label{mult by w2/psi}
\begin{array}{ll}
\ds{\frac{\mu(\lambda',M)}{\lambda'}\int_{C}w^2}=&\underbrace{\frac{1}{\lambda'}\int_{C}\left|\frac{\nabla\psi^{\lambda',M}}{\psi^{\lambda',M}}\,w-\nabla w+\frac{\lambda'}{M}\tilde{e}\,w\right|^2}_{\geq0}+\vspace{3 pt}\\
&\ds{\int_C (q\cdot\tilde{e})\,w^2+\frac{1}{\lambda'}\int_{C}\left[\zeta w^2-|\nabla w|^2\right]},
\end{array}
\end{equation}
for all $\lambda'>0$ and $M>0$ (in the case of a general diffusion matrix, see Remark \ref{for general diffusion A}). Since (\ref{mult by w2/psi}) is true for any $w\in \mathcal{I},$ then
$$\forall\,\lambda',M>0,~\frac{\mu(\lambda',M)}{\lambda'}\geq h(\lambda')\geq \inf_{\lambda'>0}h(\lambda').$$
Having (\ref{var form in mu}), one then concludes that for any $\lambda',M\in (0,+\infty)$ and for any $w\in \mathcal{I}$ with $||w||_{L^2(C)}=1,$
\begin{equation}\label{important ineq}
\ds{\inf_{\lambda'>0}h(\lambda')\leq\frac{c^*(M)}{M}\leq }\ds{\frac{1}{\lambda'}\int_{C}\left|\frac{\nabla\psi^{\lambda',M}}{\psi^{\lambda',M}}\,w-\nabla w+\frac{\lambda'}{M}\tilde{e}\,w\right|^2+h(\lambda').}
\end{equation}
To complete the proof we need the following
\begin{lemma}\label{convergence} Let $\{M_n\}_{n\in\N}$ be a sequence of positive real numbers such that $M_n\rightarrow+\infty$ as $n\rightarrow+\infty.$ Then, for a fixed $\lambda'>0,$ the sequence $\{\psi^{\lambda',M_n}\}_{n\in\N}$ of principal eigenfunctions
of the problem (\ref{Leq mu}) corresponding to $M=M_n$ converges, in $H^{1}_{loc}(\Omega)$ strong, to a function
$\psi^{\lambda',+\infty}\in H^{1}_{loc}(\Omega)$ as $n\rightarrow+\infty.$ Moreover, $\psi^{\lambda',+\infty}$ is a first integral of $q$
and
\begin{equation}\label{D}
\ds{\lim_{n\rightarrow+\infty}\int_C\left|\frac{\nabla \psi^{\lambda',M_n}}{\psi^{\lambda',M_n}}\psi^{\lambda',+\infty}-\nabla\psi^{\lambda',+\infty}\right|^2=0.}
\end{equation}
\end{lemma}
The proof of this lemma will be postponed at the end of the proof of Theorem \ref{main}.

Now, we consider any sequence $\{M_n\}_{n\in\N}$ in $(0,+\infty)$ such that $M_n\rightarrow+\infty$ as $n\rightarrow+\infty.$
 Before going further on, we mention that parts (iv) and (v) of Proposition \ref{prop} yield that
\begin{equation}\label{inf h = max on I1}
\ds{\max_{w\in \mathcal{I}_1}\frac{\int_{C}q\cdot\tilde{e}\,w^2}{\int_Cw^2}=\inf_{\lambda'>0} h(\lambda')}.
\end{equation}
Together with (\ref{important ineq}), we consequently have
\begin{equation}\label{liminf c*M}
\ds{\max_{w\in \mathcal{I}_1}\frac{\int_{C}q\cdot\tilde{e}\,w^2}{\int_Cw^2}\leq\liminf_{n\rightarrow+\infty}\frac{c^{*}(M_n)}{M_n}.}
\end{equation}

Part (iii) of Proposition \ref{prop} and (\ref{important ineq}) lead us to two different cases according to the nature of the function $h.$ The first case is when $h$ is convex and decreasing on $(0,+\infty).$ We apply (\ref{important ineq}) for $w=\psi^{\lambda',+\infty}$ (where $\lambda'>0$ is
arbitrarily chosen) together with (\ref{D}) and we get
$$\begin{array}{ll}
\ds{\limsup_{n\rightarrow+\infty} \frac{ c^*(M_n)}{M_n}\leq}&\ds{\frac{1}{\lambda'}\lim_{n\rightarrow+\infty}\int_{C}\left|\frac{\nabla\psi^{\lambda',M_n}}{\psi^{\lambda',M_n}}\,\psi^{\lambda',+\infty}-\nabla \psi^{\lambda',+\infty}\right|^2}+h(\lambda')=h(\lambda')
\end{array}$$
by (\ref{D}). Since this is true for any $\lambda'>0,$ then
\begin{equation}\label{limsup c*M}
\ds{\limsup_{n\rightarrow+\infty} \frac{ c^*(M_n)}{M_n}\leq \lim_{\lambda\rightarrow+\infty}h(\lambda)=\max_{w\in \mathcal{I}_1}\frac{\int_{C}q\cdot\tilde{e}\,w^2}{\int_Cw^2}}
\end{equation}
by part (iv) of the proposition. From (\ref{liminf c*M}) and (\ref{limsup c*M}), we get the result in the first case.

The second case is when the function $h$ attains its minimum at $\lambda_0>0.$ We apply (\ref{important ineq}) and (\ref{D}) for $\lambda'=\lambda_0,$
$w=\psi^{\lambda_0,+\infty},$ and $M=M_n.$ Hence,
$$\ds{\limsup_{n\rightarrow+\infty} \frac{ c^*(M_n)}{M_n}=h(\lambda_0)=\min_{\lambda'>0}h(\lambda')}.$$
Part (v) of Proposition \ref{prop} together with (\ref{liminf c*M}) (which is true in both cases) yield that
$$\ds{\lim_{n\rightarrow+\infty} \frac{ c^*(M_n)}{M_n}=\max_{w\in \mathcal{I}_1}\frac{\int_{C}q\cdot\tilde{e}\,w^2}{\int_Cw^2}},$$
in the second case.

Thus, in both cases, the limit of $c^*(M_n)/{M_n}$ is the same. Moreover, this limit is obtained for an arbitrarily
chosen sequence $\{M_n\}_n$ converging to $+\infty$ as $n\rightarrow+\infty.$ This implies that
$\lim_{M\rightarrow+\infty}c^{*}(M)/M$ exists and it is equal to $$\ds{\max_{w\in \mathcal{I}_1}\frac{\int_{C}q\cdot\tilde{e}\,w^2}{\int_Cw^2}},$$
which eventually proves Theorem \ref{main}.\hfill$\Box$

\vskip0.5cm
Now, we turn to prove Lemma \ref{convergence} which was announced and used in the proof of Theorem \ref{main}.
\vskip0.5cm
\textbf{Proof of Lemma \ref{convergence}.} Let us fix $\lambda'>0$ and take any sequence of positive real numbers $\{M_n\}_{n\in\N}$ converging to $+\infty$ as $n\rightarrow+\infty$. For any $n\in\N,$ the principal eigenfunction $\psi^{\lambda',M_n}$
is unique up to multiplication by a nonzero constant. Hence, we can assume that
\begin{equation}\label{normalization}
\forall\,n\in\N,~~\int_C\left(\psi^{\lambda',M_n}\right)^2=1.
\end{equation}
 We multiply (\ref{Leq mu}) (where $M=M_n$) by
$\psi^{\lambda',M_n}$ and we integrate by parts over the periodicity cell $C.$ Owing to the periodicity of $q$ and $\zeta$ together with
the condition (\ref{cq}), we then get
\begin{equation}\label{mult by psi}
-\int_{C}|\nabla\psi^{\lambda',M_n}|^2+\left(\frac{\lambda'}{M_n}\right)^2+\int_{C}\left[\lambda'q\cdot\tilde{e}+\zeta\right]\left(\psi^{\lambda',M_n}\right)^2
=\mu(\lambda',M_n)
\end{equation}
for all $n\in\N.$ As direct consequences of (\ref{mult by psi}), we have (from a certain rank $n_0$)
\begin{equation}\label{mu is bounded}
\forall n\geq n_0,~0<\mu(\lambda',M_n)\leq \lambda'^2+\lambda'\|q\cdot\tilde{e}\|_{\infty}\,+\|\zeta\|_{\infty},
\end{equation}
and $$\forall n\geq n_0,~0<\frac{\mu(\lambda',M_n)}{\lambda'}\leq \lambda'+\|q\cdot\tilde{e}\|_{\infty}+\frac{\|\zeta\|_{\infty}}{\lambda'}.$$
In other words, the sequences $\left\{\mu(\lambda',M_n)\right\}_{n\in\N}$ and $\left\{\frac{\mu(\lambda',M_n)}{\lambda'}\right\}_{n\in\N}$ are bounded whenever $\lambda'>0$ is fixed. Thus there exists $\mu(\lambda',+\infty)\geq0$ such that, up to extraction of a subsequence,
\begin{equation}
\mu(\lambda',M_n)\rightarrow\mu(\lambda',+\infty) \hbox{ as }n\rightarrow+\infty.
\end{equation}
Equations (\ref{normalization}) and (\ref{mult by psi}) imply that the sequence $\left\{\psi^{\lambda',M_n}\right\}_{n\in\N}$ is bounded in $H^1(C).$
It follows that there exists $\psi^{\lambda',+\infty}\in H^{1}_{loc}(\Omega)$ such that, up to extraction of a subsequence, $\psi^{\lambda',M_n}\rightarrow\psi^{\lambda',+\infty}$ in
$H^1_{loc}(\Omega)$ weak, in $L^2_{loc}(\Omega)$ strong, and almost everywhere in $\Omega$ as $n\rightarrow+\infty.$ Thus, the function $\psi^{\lambda',+\infty}$ is $L$-periodic with respect to $x.$ The strong convergence in $L^2(C)$ leads to $$\int_C\left(\psi^{\lambda',+\infty}\right)^2=1$$ and hence $\psi^{\lambda',+\infty}\not\equiv 0$ in $\Omega.$
We replace $M$ by $M_n$ in (\ref{Leq mu}), divide the equation by $M_n$ and we pass to the limit as $n\rightarrow+\infty$ in the sense of distributions. From the weak convergence of $\left\{\psi^{\lambda',M_n}\right\}_{n\in\N}$ and the boundedness of $\{\mu(\lambda',M_n)\}_{n\in\N},$ one then has $q\cdot \nabla \psi^{\lambda',+\infty}=0$ in $\mathcal{D}'(\Omega)$. Consequently, $q\cdot \nabla \psi^{\lambda',+\infty}=0$ almost everywhere in $\Omega.$ That is, $\psi^{\lambda',+\infty}$ is a nonzero first integral of $q.$

Now, we multiply (\ref{Leq mu}) (where $M=M_n$) by $\psi^{\lambda',+\infty},$ we integrate by parts over $C$ and we pass to the limit as $n\rightarrow+\infty.$ We notice that
$$\ds{\int_{C}\Delta\psi^{\lambda',M_n}\psi^{\lambda',+\infty}}=\ds{-\int_{C}\nabla\psi^{\lambda',M_n}\cdot\nabla\psi^{\lambda',+\infty}-\displaystyle{\frac{\lambda'}{M_n}\int_{\partial C}\nu\cdot \tilde{e}\,\psi^{\lambda',M_n}\psi^{\lambda',+\infty}}\rightarrow-\int_{C}\left|\nabla\psi^{\lambda',+\infty}\right|^2}$$ \hbox{ as $n\rightarrow+\infty$} (from the strong convergence in $L^2_{loc}(\Omega)$, the boundary term converges to $0$. We use the weak convergence in $H^1_{loc}(\Omega)$ for the limit of the first term), and
$$
    \begin{array}{rl}
\ds{\int_{C}\left(q\cdot\nabla\psi^{\lambda',M_n}\right)\psi^{\lambda',+\infty}}=&\ds{-\int_{C}\left(\nabla\cdot q\right)\,\psi^{\lambda',+\infty}\psi^{\lambda',M_n}-\int_{C}\left(q\cdot\nabla \psi^{\lambda',+\infty}\right)\psi^{\lambda',M_n}}\vspace{3pt}\\
=&0\hbox{ (from (\ref{cq}) and since $\psi^{\lambda',+\infty}\in\,\mathcal{I}$).}
    \end{array}
$$
Hence,
\begin{equation}\label{mult by psi infty and pass to limit}
\mu(\lambda',+\infty)=-\int_{C}\left|\nabla\psi^{\lambda',+\infty}\right|^2+\int_{C}\left[\lambda'q\cdot\tilde{e}+\zeta\right]\left(\psi^{\lambda',+\infty}\right)^2.
\end{equation}
On the other hand, we take $w=\psi^{\lambda',+\infty}$ and $M=M_n$ in (\ref{mult by w2/psi}), we multiply the equation by the fixed $\lambda'>0$
and we pass to the limit as $n\rightarrow+\infty$ to obtain
\begin{equation}
\begin{array}{c}
D:=\ds{\lim_{n\rightarrow+\infty}\int_{C}\left|\frac{\nabla\psi^{\lambda',M_n}}{\psi^{\lambda',M_n}}\,\psi^{\lambda',+\infty}-\nabla \psi^{\lambda',+\infty}+\frac{\lambda'}{M_n}\tilde{e}\,\psi^{\lambda',+\infty}\right|^2}\vspace{3pt}\\
\ds{=\mu(\lambda',+\infty)-\lambda'\int_C q\cdot\tilde{e}\,{\left(\psi^{\lambda',+\infty}\right)}^2-\int_{C}\left[\zeta {\psi^{\lambda',+\infty}}^2-|\nabla \psi^{\lambda',+\infty}|^2\right].}
\end{array}
\end{equation}
 However, $\ds{\frac{\lambda'}{M_n}\tilde{e}\,\psi^{\lambda',+\infty}\rightarrow0}$ in $L^2(C)$ strong as $n\rightarrow+\infty.$ Also,
$\ds{\left\{\left|{\nabla\psi^{\lambda',M_n}}/{\psi^{\lambda',M_n}}\right|\right\}_{n\in\N}}$ is bounded in $L^2(C)$ (we simply divide (\ref{Leq mu}) by $\psi^{\lambda',M_n}$ and integrate by parts over $C$. This leads to $$\ds{\int_{C}\left|\frac{\nabla\psi^{\lambda',M_n}}{\psi^{\lambda',M_n}}\right|^2+\left(\frac{\lambda'}{M_n}\right)^2|C|+\int_{C}\zeta}=\mu(\lambda',M_n)|C|;$$
hence, $\ds{\left\{\int_{C}\left|\frac{\nabla\psi^{\lambda',M_n}}{\psi^{\lambda',M_n}}\right|^2\right\}_n}$ is bounded due to (\ref{mu is bounded})).
Consequently,
\begin{equation}\label{D=}
\ds{D=\lim_{n\rightarrow+\infty}\int_{C}\left|\frac{\nabla\psi^{\lambda',M_n}}{\psi^{\lambda',M_n}}\,\psi^{\lambda',+\infty}-\nabla \psi^{\lambda',+\infty}\right|^2.}
\end{equation}
 Referring to (\ref{mult by psi infty and pass to limit}), we finally obtain
$$\lim_{n\rightarrow+\infty}\int_{C}\left|\frac{\nabla\psi^{\lambda',M_n}}{\psi^{\lambda',M_n}}\,\psi^{\lambda',+\infty}-\nabla \psi^{\lambda',+\infty}\right|=0.$$
Since $\psi^{\lambda',M_n}\rightarrow\psi^{\lambda',+\infty}$ in $L^2_{loc}(\Omega)$ strong, (\ref{mult by psi}) yields that $\ds{\int_{C}|\nabla\psi^{\lambda',M_n}|^2}$ converges to $$\ds{\lambda'\int_{C}q\cdot\tilde{e}\left(\psi^{\lambda',+\infty}\right)^2+\int_{C}\zeta\left(\psi^{\lambda',+\infty}\right)^2-\mu(\lambda',+\infty)}
.$$ Equation (\ref{mult by psi infty and pass to limit}) again yields that
\begin{equation}\label{strong conv in H1 with A}
\ds{\int_{C}|\nabla\psi^{\lambda',M_n}|^2\rightarrow\ds{\int_{C}|\nabla\psi^{\lambda',+\infty}|^2}}\hbox{ as }n\rightarrow+\infty.
\end{equation}
Eventually, $\{\psi^{\lambda',M_n}\}_n$ converges to $\psi^{\lambda',+\infty}$ in $H^1_{loc}(\Omega)$ strong and this completes the proof of Lemma \ref{convergence}.\hfill$\Box$

\begin{remark}[about the proof of (\ref{large advection}) in case of a  diffusion  $A=A(x,y)$]\label{for general diffusion A}
\rm{ In this remark, we detail some differences
which arise in the proof of (\ref{large advection}) when we consider, instead of the identity matrix, a general diffusion matrix $A(x,y)$ satisfying (\ref{cA}). The eigenvalue problem (\ref{Leq mu}) becomes
\begin{equation}\label{Leq mu with A}
\left\{
  \begin{array}{rl}
 \mu(\lambda',M)\psi^{\lambda',M}   =&\displaystyle{\nabla\cdot A\nabla \psi^{\lambda',M}+2\frac{\lambda'}{M}\tilde{e}\cdot
A\nabla\psi^{\lambda',M}+M\,q\cdot\nabla\psi^{\lambda',M}}\vspace{3pt}\\
&+\displaystyle{\left[\left(\frac{\lambda'}{M}\right)^2\tilde{e}\cdot A\tilde{e}+\lambda'\,
q\cdot\tilde{e}+\frac{\lambda'}{M}\nabla\cdot A\tilde{e}+\zeta\right]\psi^{\lambda',M}}\hbox{ in }\Omega,\vspace{3 pt}\\
\nu\cdot A\nabla\psi^{\lambda',M}=&\ds{-\frac{\lambda'}{M}(\nu \cdot
A\tilde{e})\psi^{\lambda',M}}\;\hbox{on}\;\partial{\Omega}\hbox{ (whenever }\partial\Omega\neq\emptyset\hbox{)}.
  \end{array}
\right.
\end{equation}
Similar to the case where $A=Id,$ the principal eigenfunctions of (\ref{Leq mu with A}) $\psi^{\lambda',M}$ will be positive in $\overline{\Omega}$, unique up to multiplication by
a nonzero constant, and $L$-periodic with respect to $x$ (see section 5 in \cite{bh}). Consequently, for any sequence $\{M_n\}_n$ in $(0,+\infty)$ such that $M_n\rightarrow+\infty$ as $n\rightarrow+\infty,$ multiplying (resp. dividing) (\ref{Leq mu with A}), for any $\lambda'>0$ and for $M=M_n,$ by $\psi^{\lambda',M_n}$ and integrating by parts over the cell $C$ implies the boundedness of $\{\nabla\psi^{\lambda',M_n}\}_n$( resp. $\{\nabla\psi^{\lambda',M_n}/\psi^{\lambda',M_n}\}_n$) in  $L^2_{loc}(\Omega).$ On the other hand, the equation (\ref{mult by w2/psi}), which was essential in the proof when $A=Id,$ will be replaced by
\begin{equation}\label{mult by w2/psi with A}
\begin{array}{ll}
\ds{\frac{\mu(\lambda',M)}{\lambda'}\int_{C}w^2}=&\ds{\frac{D_A(\lambda',M)}{\lambda'}}
\ds{+\int_C (q\cdot\tilde{e})\,w^2+\frac{1}{\lambda'}\int_{C}\left[\zeta w^2-\nabla w\cdot A\nabla w\right],}
\end{array}
\end{equation}
for all $\lambda'$ and $M$ in $(0,+\infty)$ where $$D_A(\lambda',M):=\int_{C}\left(\frac{\nabla\psi^{\lambda',M}}{\psi^{\lambda',M}}\,w-\nabla w+\frac{\lambda'}{M}\tilde{e}\,w\right)\cdot A\left(\frac{\nabla\psi^{\lambda',M}}{\psi^{\lambda',M}}\,w-\nabla w+\frac{\lambda'}{M}\tilde{e}\,w\right).$$ The result of Lemma \ref{convergence} will remain true with (\ref{D}) replaced by
$$\lim_{n\rightarrow+\infty}D_A(\lambda',M_n)=0.$$
Finally, due to the coercivity of the diffusion matrix $A$ given in (\ref{cA}), we can be easily adapt the proof of the results of Proposition \ref{prop} to the
function $\lambda\mapsto h^A(\lambda)$ defined by \begin{equation}\label{hA}
\forall\lambda\in(0,+\infty),~\ds{h^A(\lambda):=\int_C (q\cdot\tilde{e})\,w^2+\frac{1}{\lambda}\int_{C}\left[\zeta w^2-\nabla w\cdot A\nabla w\right],}
\end{equation}
which coincides with $\lambda\mapsto h(\lambda)$ defined in (\ref{defh}) when $A=Id.$
}
\end{remark}

\subsection{Cases of large advection with small reaction or large diffusion (proof of Theorem \ref{mixing large advection with diff and reac})}\label{proof in mixed cases}

We mention that the proof of (\ref{large advection large diffusion}) is very similar to that of (\ref{large advection small reaction}). We are going to prove the limit in (\ref{large advection small reaction}) only.\vskip 0.35cm

\textbf{\underline{Step 1. Existence of a maximizer for (\ref{large advection small reaction})}:}
At the beginning, we prove that  $\ds\sup_{\ds{w\in \mathcal {I}}}\frac{ \ds\int_C(q\cdot\tilde{e})\,w}{\ds\sqrt{\int_C \nabla w\cdot A\nabla w}}$ is finite. For any $w\in\mathcal{I},$ we define $$\overline{w}:=\mi_{\!\!\!\!C}w(x)dx,$$ and we write $w=\overline{w}+v.$ We notice that $\nabla w\equiv\nabla v$ and, thanks to Poincar\'e's inequality,
we get
\begin{equation}\label{poincare}
\ds{||v||_{_{L^2(C)}}^2}\leq \kappa\int_C |\nabla v|^2=\kappa\int_C |\nabla w|^2\leq \frac{\kappa}{\alpha_1} \int_C \nabla w\cdot A\nabla w,
 \end{equation}
for some $\kappa>0$ independent of $w$ and $v,$ where $\alpha_1>0$ is given by (\ref{cA}). Moreover, it follows from the fourth line in (\ref{cq}) that $$\int_Cq\cdot\tilde{e}\,w=\int_Cq\cdot\tilde{e}\,v.$$
Thus, applying Cauchy-Schwarz inequality, we get
$$\forall w\in\mathcal{I},~~\left|\int_Cq\cdot\tilde{e}\,w\right|=\left|\int_Cq\cdot\tilde{e}\,v\right|\leq||q\cdot\tilde{e}||_{_{L^2(C)}}||v||_{_{L^2(C)}}\leq\sqrt{\kappa/\alpha_1}||q\cdot\tilde{e}||_{_{L^2(C)}}\sqrt{\int_C \nabla w\cdot A\nabla w}.$$
Hence, for any $w\in\mathcal{I}$ such that $w$ is not constant,
$$0\leq\frac{\left|\int_Cq\cdot\tilde{e}\,w\right|}{\sqrt{\int_C \nabla w\cdot A\nabla w}}\leq \sqrt{\frac{\kappa}{\alpha_1}}||q\cdot\tilde{e}||_{_{L^2(C)}}<+\infty$$ since $q\in C^{1,\delta}(\overline{\Omega}).$ Consequently, the quantity $$l:=\ds\sup_{\ds{w\in \mathcal {I}}}\frac{ \ds\int_C(q\cdot\tilde{e})\,w}{\ds\sqrt{\int_C \nabla w\cdot A\nabla w}}\geq0$$ is well defined.
We mention that if $\int_C(q\cdot\tilde{e})\,w=0$ for each $w\in\mathcal {I},$ then $l=0$ and the supremum is a maximum attained by any nonconstant $w\in\mathcal {I}.$ In what follows, we have to treat the case where there exists at least a  $w_0\in \mathcal {I}$ (nonconstant a.e in $C$) such that $\int_C(q\cdot\tilde{e})\,w_0\neq0$ and consequently $l>0.$

Now, we prove that the above supremum is actually a maximum. We take $\{w_n\}_{n\in\N}$ as a maximizing sequence. As it was done in (\ref{poincare}) above, we may assume that
\begin{equation}\label{poincare 1}
\ds{\forall\,n\in\N,~~ ||w_n||_{_{L^2(C)}}^2}\leq \kappa\int_C |\nabla w_n|^2\leq \frac{\kappa}{\alpha_1} \int_C \nabla w_n\cdot A\nabla w_n.
 \end{equation}
Moreover, we can consider $$\tilde {w}_n:=\frac{w_n}{\sqrt{\int_C \nabla w_n\cdot A\nabla w_n}}$$
as a maximizing sequence. The advantage is that $\{\tilde {w}_n\}_n$ is bounded in $L^2(C)$ (due to (\ref{poincare 1})) and in $H^1(C)$ since
$\int_C \nabla \tilde w_n\cdot A\nabla \tilde w_n=1$ for each $n\in\N.$ As a consequence, there exists $\tilde w\in H^1(C)$
such that
$$\tilde w_n\rightarrow\tilde w~\hbox{ in $L^2(C)$ strong and }~\tilde w_n\rightharpoonup\tilde w ~\hbox{ in $H^1(C)$ weak as $n\rightarrow+\infty$.}$$
Thus, \begin{equation}\label{liminf weak}
\ds1=\liminf_{n\rightarrow+\infty}\int_C \nabla \tilde w_n\cdot A\nabla \tilde w_n\geq \int_C \nabla \tilde w\cdot A\nabla \tilde w,
\end{equation}
and $\tilde{w}\in\mathcal{I}$ as the weak limit of first integrals of $q.$\\
On the other hand, strong convergence in $L^2(C),$  the definition of the maximizing sequence $\{\tilde w_n\}_n$ and (\ref{liminf weak}) yield that
\begin{equation}\label{moh}
\begin{array}{lll}
  \ds\int_C(q\cdot \tilde e)\tilde w &=&\ds\lim_{n\rightarrow+\infty} \int_C(q\cdot \tilde e)\tilde w_n=l\geq l\times \int_C \nabla \tilde w\cdot A\nabla \tilde w.
\end{array}
\end{equation}
We mention that $\tilde w$ cannot be constant almost everywhere in $C$ because in this case one gets $l=\int_C(q\cdot \tilde e)\tilde w=\tilde w\int_Cq\cdot\tilde e =0,$ and this contradicts the assumption that $l>0.$\\
Therefore, it follows from (\ref{moh}) and the definition of $l,$ that $$\ds{l=\frac{\int_C(q\cdot \tilde e)\tilde w }{\sqrt{\int_C \nabla \tilde w\cdot A\nabla \tilde w}}},$$ and so, the maximum of (\ref{large advection small reaction}) is attained at $\tilde w.$

We also obtain $\int_C \nabla \tilde w\cdot A\nabla \tilde w=1$ which yields that $\{\tilde w_n\}_n$ converges to $\tilde w$ in $H^1(C)$
\textbf{strong} as $n\rightarrow+\infty$.
\vskip 0.5cm
\underline{\textbf{Step 2.}}
Theorem \ref{main} yields that, for any $\epsilon>0,$ the limit of ${\ds{c^{*}_{\Omega,A,M\,q,\epsilon f}(e)}}/{M}$
is related to the set $$\mathcal{I}_1^{A,\epsilon}:=\left\{w\in\,\mathcal{I},\hbox{ such that }\ds\epsilon\!\!\int_C\zeta w^2\geq\int_C\nabla w\cdot A\nabla w\right\}.$$

As we did at the beginning of Step 1, we write each $w\in \mathcal{I}_1^{A,\epsilon}$ as $w=v+\overline {w},$ where $$\overline {w}=\ds\frac{\ds\int_Cw}{|C|}.$$
For each  $w\in \mathcal{I}_1^{A,\epsilon}$ with $||w||_{L^2(C)}=1,$ we have  $|\overline{w}|\leq1/\sqrt{|C|}.$ We also have $|\int_{C}v|\leq\sqrt{|C|}||v||_{L^2(C)}.$ Owing to (\ref{poincare})  together with the facts that $w\in\mathcal{I}_1^{A,\epsilon}$ and $\zeta$ is globally bounded, one consequently gets
\begin{equation}\label{v is a big O}
\int_{C}v=O(\sqrt{\epsilon})\hbox{ as }\epsilon\rightarrow0^+.
\end{equation}
On the other hand, $$1=\int_Cw^2=\overline{w}^2|C|+2\underbrace{\overline{w}\int_{C}v}_{O(\sqrt{\epsilon})}+\underbrace{\int_Cv^2}_{O(\epsilon)}.$$
Thus, $\overline{w}^2|C|=1+O(\sqrt{\epsilon})$ as $\epsilon\rightarrow0^+.$ Now, we write $$\overline{w}^2|C|-1=(\overline{w}\sqrt{|C|}-1)(\overline{w}\sqrt{|C|}+1)=O(\sqrt{\epsilon}),$$
and we use the fact that $$\begin{array}{rl}
1\leq\overline{w}\sqrt{|C|}+1\leq 2&\hbox{ when } \overline{w}\geq0,\hbox{ and}\vspace{4 pt}\\
-2\leq\overline{w}\sqrt{|C|}-1\leq -1 &\hbox{ when } \overline{w}\leq0,
\end{array}$$ to obtain
\begin{equation}\label{w bar is 1+ a big O}
\forall \,w\in \mathcal{I}_1^{A,\epsilon},~~ (\;||w||_{L^2(C)}=1\;)\Rightarrow\;\overline{w}=\text{sgn}(\overline{w})\frac{1}{\sqrt{|C|}}+O(\sqrt{\epsilon})~\hbox{ as }~\epsilon\rightarrow0^+.
\end{equation}
For such $w$'s, having $q\cdot\tilde{e}\in L^\infty(C),$ it then follows from (\ref{v is a big O}) and (\ref{w bar is 1+ a big O}) that $$\begin{array}{ll}
\ds\int_{C}(q\cdot\tilde{e})w^2=2\,\overline{w}\ds\int_C(q\cdot \tilde{e})v+\ds\int_C(q\cdot \tilde{e})v^2&=\ds\frac{2\text{ sgn}(\overline{w})\ds\int_C(q\cdot \tilde{e})v\,dx}{\ds\sqrt{|C|}}+O(\epsilon)~\hbox{ as }~\epsilon\rightarrow0^+\vspace{4pt}\\
&=\ds\frac{2\text{ sgn}(\overline{w})\ds\int_C(q\cdot \tilde{e})w\,dx}{\ds\sqrt{|C|}}+O(\epsilon),
\end{array}$$
since $\ds\int_C(q\cdot \tilde{e})=0.$\\
As $w\in \mathcal{I}_1^{A,\epsilon},$ we have $\ds{\frac{\sqrt{\epsilon}\sqrt{\int_C\zeta w^2}}{\sqrt{\int_C\nabla w\cdot A\nabla w}}}\geq 1.$ Hence,
\begin{equation}\label{key ineq for mixed drift and reaction}
\forall\, w\in \mathcal{I}_1^{A,\epsilon},\;\ds\int_{C}(q\cdot\tilde{e})w^2\leq \ds{\frac{\sqrt{\epsilon}\sqrt{\int_C\zeta w^2}}{\sqrt{\int_C\nabla w\cdot A\nabla w}}}\times\ds\frac{2\left|\int_C(q\cdot \tilde{e})w\,dx\right|}{\ds\sqrt{|C|}}+O(\epsilon),
\end{equation}
whenever $||w||_{L^2(C)}=1.$

Moreover, if $\int_C\nabla w\cdot A\nabla w={\epsilon\!\!\int_C\zeta\,w^2}~ \hbox{ and }~||w||_{L^2(C)}=1,$ then we have:
\begin{equation}\label{key equal for mixed drift and reaction}
\begin{array}{ll}
\ds\int_{C}(q\cdot\tilde{e})w^2= \ds{\frac{\sqrt{\epsilon}\sqrt{\int_C\zeta w^2}}{\sqrt{\int_C\nabla w\cdot A\nabla w}}}\times\ds\frac{2\text{ sgn}(\overline{w})\int_C(q\cdot \tilde{e})w\,dx}{\ds{\sqrt{|C|}}}+O(\epsilon).
\end{array}
\end{equation}
\vskip 0.5cm

\underline{\textbf{Step 3.}} For a fixed $\epsilon>0,$ we know from Theorem \ref{main} that
\begin{equation}\label{first}
\lim_{M\rightarrow+\infty}\ds{\frac{\ds{c^{*}_{\Omega,A,M\,q,\epsilon f}(e)}}{M\sqrt{\epsilon}}}=\max_{\hskip0.7cm{w\in \mathcal{I}_1^{A,\epsilon},}_{_{_{\hskip-1.5cm \ds||w||_{L^2(C)}=1}}}}\ds\frac{1}{\sqrt{\epsilon}}\int_{C}(q\cdot\tilde{e})w^2=\frac{1}{\sqrt{\epsilon}}\int_{C}(q\cdot\tilde{e})w_\epsilon^2,
\end{equation}
for some $w_{\epsilon}\in \mathcal{I}_1^{A,\epsilon}$ with $||w_{\epsilon}||_{L^2(C)}=1.$ However,
$$||\nabla w_\epsilon||_{L^2(C)}^2\leq\frac{1}{\alpha_1}\int_C\nabla w_\epsilon\cdot A\nabla w_\epsilon\leq\frac{\epsilon||\zeta||_\infty}{\alpha_1},$$
since $w_{\epsilon}\in \mathcal{I}_1^{A,\epsilon}.$ Hence, $\nabla w_\epsilon\rightarrow 0$ in $L^2(C)$ as $\epsilon\rightarrow0^+.$
Consequently, the sequence $\{w_\epsilon\}_{\epsilon}$ is bounded in $H^1(C),$ and thus, there exists $w_0\in H^1(C)$ such that $$w_\epsilon\rightharpoonup w_0\hbox{ in }H^1(C) \hbox{ weak, and }w_\epsilon\rightarrow w_0 \hbox{ in }L^2(C) \hbox{ strong, as }\epsilon\rightarrow 0^+.$$
Strong convergence in $L^2(C)$ yields that $w_0\not\equiv0$ in $C$ and $||w_0||_{L^2(C)}=1.$ Besides, weak convergence implies that $\int_C|\nabla w_0|^2\leq\liminf_{\epsilon\rightarrow0^+}\int_C|\nabla w_\epsilon|^2=0.$ Therefore, $|\nabla w_0|=0$ a.e in $C$ and
$\ds{w_0=\frac{\pm1}{\sqrt{|C|}}}$ is constant a.e in $C.$ One concludes that $\int_C\zeta w_\epsilon^2\rightarrow\frac{1}{|C|}\int_C\zeta$
as $\epsilon \rightarrow 0^+.$ These results together with (\ref{key ineq for mixed drift and reaction}) and (\ref{first}) lead to
\begin{equation}\label{limsup c*eps M}
 \limsup_{\epsilon\rightarrow0^+}\lim_{M\rightarrow+\infty}\ds{\frac{\ds{c^{*}_{\Omega,A,M\,q,\epsilon f}(e)}}{M\sqrt{\epsilon}}}\leq \frac{2\sqrt{\int_C\zeta}}{|C|}\ds\max_{\ds{w\in \mathcal {I}}}\frac{ \int_C(q\cdot\tilde{e})\,w}{\sqrt{\int_C \nabla w\cdot A\nabla w}}.
\end{equation}

\underline{\textbf{Step 4.}} To prove equality, we take  any maximizer $v\in\mathcal{I}$ (non constant) of $$R(w):=\frac{ \int_C(q\cdot\tilde{e})\,w}{\sqrt{\int_C \nabla w\cdot A\nabla w}},$$ over $\mathcal{I}.$ It is easy to see that for any $k\in\R,$ $v+k$ is also
 a maximizer of $R.$ Thus, we can choose without any loss of generality  $v$ so that $\overline{v}\geq0.$

We want to prove that there exists $\epsilon_0>0,$ such that for any $0<\epsilon\leq\epsilon_0$ we can find a maximizer $w_\epsilon\in\mathcal{I}$ of $R$  satisfying
\begin{equation}\label{normalized cond}
||w_{\epsilon}||_{L^2(C)}=1\text{ and }
\ds{\int_C\nabla w_\epsilon\cdot A\nabla w_\epsilon={\epsilon\!\!\int_C\zeta\,w_\epsilon^2}.}
\end{equation}
Using Poincar\'e's inequality (as in Step 1), we have $\int_C(v-\overline{v})^2\leq\kappa\int_C|\nabla v|^2$ for some $\kappa>0$ depending
only on the set $C.$ Moreover, the function $\zeta$ is positive and belongs to $L^\infty(C).$ Thus, for any $0<\epsilon\leq\epsilon_0:=\frac{\alpha_1}{\kappa||\zeta||_\infty}$ (where $\alpha_1$ is given by (\ref{cA})), we have
\begin{equation}\label{normali 1}
\epsilon\int_C\zeta(v-\overline{v})^2\leq\int_C\nabla v\cdot A\nabla v.
\end{equation}
 Now, let's fix $\epsilon$ in $(0,\epsilon_0].$ Since $\overline{v}\geq0,$ one can then find
a constant $m=m(\epsilon,v)\leq0$ depending on $\epsilon$ and $v$ such that
\begin{equation}\label{normali 2}
\epsilon\int_C\zeta(v-m)^2\geq\epsilon\left(\min_{z\in C}\zeta(z)\right)\int_C(v-m)^2\geq\int_C\nabla v\cdot A\nabla v.
\end{equation}
The continuity of $t\mapsto \epsilon\int_C(v-t)^2-\int_C\nabla v\cdot A\nabla v$ together with (\ref{normali 1}) and (\ref{normali 2}) yield that there exists $r=r(\epsilon,v)\in[m,\overline{v}],$ such that $$\epsilon\int_C\zeta(v-r)^2=\int_C\nabla v\cdot A\nabla v.$$
We call $$\ds{w_\epsilon:=\ds{\frac{v-r}{||v-r||_{L^2(C)}}}}.$$
Then, $w_\epsilon$ verifies (\ref{normalized cond}) and it maximizes $R(w)$ since $v$ does and since $\ds{\int_C q\cdot\tilde{e}=0}.$ Moreover,
$\overline{w}_\epsilon\geq0$ since $r\in[m,\overline{v}].$
 Imitating the argument used in Step 3, one gets that, up to the extraction of a subsequence, $w_\epsilon\rightarrow \frac{1}{\sqrt{|C|}}$ strongly in $L^2(C).$

Applying (\ref{key equal for mixed drift and reaction}) for $w=w_\epsilon$ (where $\text{sgn}(\overline{w}_\epsilon)=+1$) and since $w_\epsilon\in\mathcal{I}^{A,\epsilon}_1,$
 we then get for any $0<\epsilon\leq\epsilon_0,$
$$\begin{array}{ll}
\ds{\max_{w\in\mathcal{I}}R(w)}=R(w_\epsilon)&=\ds\frac{\sqrt{|C|}}{2\sqrt{\epsilon\int_C\zeta w_\epsilon^2}}\int_{C}(q\cdot\tilde{e})w_\epsilon^2+O(\sqrt\epsilon)\\
&\leq O(\sqrt\epsilon)+\ds\frac{\sqrt{|C|}}{2\sqrt{\int_C\zeta w_\epsilon^2}} \ds{\max_{\hskip0.7cm{w\in \mathcal{I}_1^{A,\epsilon},}_{_{_{\hskip-1.5cm \ds||w||_{L^2(C)}=1}}}}\ds\frac{\int_{C}(q\cdot\tilde{e})w^2}{\sqrt{\epsilon}}}.
\end{array}$$
In other words,
\begin{equation}\label{key}
\forall \,0<\epsilon\leq\epsilon_0,~~\ds{\max_{w\in\mathcal{I}}R(w)}\leq O(\sqrt\epsilon)+\ds\frac{\sqrt{|C|}}{2\sqrt{\int_C\zeta w_\epsilon^2}}\lim_{M\rightarrow+\infty}\ds{\frac{\ds{c^{*}_{\Omega,A,M\,q,\epsilon f}(e)}}{M\sqrt{\epsilon}}}.
\end{equation}

Passing to the limit as $\epsilon\rightarrow0^+$ in (\ref{key}) and using the strong convergence of $w_\epsilon$ in $L^2(C),$ we  obtain
\begin{equation}\label{liminf c*eps M}
\ds{\max_{w\in\mathcal{I}}R(w)}\leq 0+\frac{|C|}{2\sqrt{\int_C\zeta}}\liminf_{\epsilon\rightarrow0^+}\lim_{M\rightarrow+\infty}\ds{\frac{\ds{c^{*}_{\Omega,A,M\,q,\epsilon f}(e)}}{M\sqrt{\epsilon}}}.
\end{equation}
This inequality and (\ref{limsup c*eps M}) finish the proof of (\ref{large advection small reaction}).\hfill$\Box$

\section{ The two dimensional case ($N=2$)}\label{N=2 section}
In this section, the space dimension is $N=2.$ In what follows, we find the form of any divergence free advection field and then prove Theorem \ref{case N=2} after passing by many auxiliary Lemmas. We first start by proving Lemma \ref{periodic trajs} which was announced in Section \ref{intro}.
\vskip0.35cm

\textbf{Proof of Lemma \ref{periodic trajs}.}
We first prove by contradiction that, for $d=1$, all the trajectories are $L_1 e_1$-periodic. Indeed, suppose that there exists a periodic unbounded trajectory $T(x)$, which is not $L_1 e_1$-periodic. Then it is $p L_1 e_1$ periodic for some $p \in \mathbb{N}$, $p\geq 2$. By periodicity of $q$, $T(x)+L_1 e_1$ is also an unbounded periodic trajectory of $q$ in $\Omega$, different from $T(x)$. Moreover, $T(x) \cap (T(x)+L_1 e_1) = \emptyset$, because two different trajectories never intersect.

We set
$$
m:= \min\left\{y_2 \text{ such that } (y_1,y_2) \in T(x) \right \},
$$
and
$$
M:= \max\left\{y_2 \text{ such that } (y_1,y_2) \in T(x) \right \}.
$$
We have $m\neq M$, otherwise $T(x)$ is a horizontal straight line and is $L_1e_1$-periodic.

We also define
$$
a_1 := \min \left\{y_1\geq0 \text{ such that } (y_1,m) \in T(x)\right\},
$$$$
a_2 :=\min \left\{y_1 \geq a_1 \text{ such that } (y_1,M) \in T(x)\right\}.
$$

Let $T'(x):= T(x)\setminus(\{(a_1,m) \}\cup\{(a_2,M)\})$. Since $T(x)$ is a simple curve, $T'(x)$ has exactly three connected components, two of which are unbounded. Let $T_b(x)$ be the bounded component of $T'(x)$, we set
$$
T_a(x) = \overline{T_b(x)},
$$
which is a compact subset of $T(x)$, with boundary $\{(a_1,m)\} \cup \{(a_2,M)\}$. We define then
$$
b_1 = \min\{y_1 \text{ such that } (y_1,y_2) \in T_a(x)\},
$$
and
$$
b_2 = \max\{y_1 \text{ such that } (y_1,y_2) \in T_a(x)\}.
$$
We define the following curve:
$$
\mathcal{C}:=T_a(x)\cup \{(a_1,y_2), \ y_2<m\}\cup\{(a_2,y_2), \ y_2>M\}.
$$
The curve $\mathcal{C}$ is a simple connected curve, which splits $\Omega$ into several connected components, two of which are unbounded. Let $\Omega_1$ be the left unbounded component, more precisely the component containing the set $\{y=(y_1,y_2) \in \Omega$ such that $y_1 < b_1\}$, and $\Omega_2$ be the right unbounded component, the one containing the set $\{y=(y_1,y_2) \in \Omega$ such that $y_1 > b_2\}$.

Let $z\in \Omega_1$ and $z' \in \Omega_2$ be two points of $T(x)+L_1e_1$, then, following $T(x)+L_1e_1$, there is a continuous path in $\Omega$ from $z$ to $z'$. This path must cross $\mathcal{C}$ for obvious reasons of continuity. However, it can not cross $\{(a_1,y_2), \ y_2<m\}$ or $\{(a_2,y_2), \ y_2>M\}$ because of the definition of $m$ and $M$.

Hence this path, which is a subset of $T(x)+L_1e_1$ must cross $T_a(x)$, which is a subset of $T(x)$. This leads to
$$
T(x) \cap (T(x)+L_1e_1) \neq \emptyset.
$$
This is a contradiction and proves the lemma for $d=1$.

\textbf{In the case} $d=2$, we need to prove that if there exists an $\mathbf{a}-$periodic unbounded trajectory, then any other unbounded periodic trajectory will be $\mathbf{a}-$periodic. The idea is to reduce this proof to the proof of the case $d=1$. Suppose that there exists an $\mathbf{a}-$periodic unbounded trajectory $T(x)$ of $q$ in $\Omega$. We set $e'_1 = \mathbf{a}/|\mathbf{a}|$, and $e'_2$ such that $(e'_1,e'_2)$ is a direct orthonormal frame of $\R^2$. In this new basis, $q$ is $L'_1-$periodic in the $e'_1$ directions, with $L'_1=|\mathbf{a}|$.

Thus, $T(x)$ is bounded in the $e'_2$ direction. Suppose now that $y=y'_1e'_1+y'_2e'_2\in \Omega$ is such that $T(y)$ is an unbounded periodic trajectory. Let $z_1 \in L_1\Z\times L_2\Z$ such that
$$
\inf_{z \in T(x)+z_1} z\cdot e'_2 > y'_2,
$$
and $z_2\in L_1\Z\times L_2\Z$ such that
$$
\sup_{z \in T(x)+z_2} z\cdot e'_2 < y'_2.
$$
We have that $T(x)+z_1$ and $T(x)+z_2$ are two trajectories of $q$ in $\Omega$, and they split $\Omega$ into three connected components, one of which is bounded in the $e'_2$ direction. We shall denote it $\Omega_b$. By construction $y \in \Omega_b$. Since two different trajectories can not intersect, $T(y)$ must stay inside $\Omega_b$, and using the same procedure as the case $d=1$ we conclude that $T(y)$ is $L'_1e'_1=\mathbf{a}-$periodic.\hfill$\Box$

\begin{proposition}\label{phodge}
Let $d=1$ or $2$ where $d$ is defined in (\ref{comega}). Let $q \in C^{1,\delta}(\overline{\Omega})$, $L$-periodic with respect to $x$ and verifying the conditions
\begin{equation} \label{q}
\left\{
\begin{array}{lll}
\displaystyle \int_C q_i & = & 0, \hbox{ for }1\leq i\leq d\leq2  \vspace{3 pt}\\
\nabla \cdot q & = & 0 \text{ in } \Omega, \vspace{3 pt}\\
q \cdot \nu & = & 0 \text{ on } \partial \Omega.
\end{array}
\right.
\end{equation}
Then, there exists $\phi \in C^{2,\delta}(\overline{\Omega})$, $L$-periodic with respect to $x$, such that
\begin{equation}\label{hodge}
q = \nabla^{\perp} \phi \ \text{ in } \Omega.
\end{equation}
Moreover, $\phi$ is constant on every connected component of $\partial \Omega$.
\end{proposition}

\begin{remark} We mention that the representation $q = \nabla^{\perp} \phi$  is already known in the case where the domain $\Omega$ is bounded and \textbf{simply connected} or  equal to whole space $\R^2.$ However, the above proposition applies in more cases due to the condition $q\cdot\nu=0$ on $\partial\Omega$ (see the proof below). For example, it applies when $\Omega$ is the whole space  $\R^2$ with a periodic array of holes or when $\Omega$ is  an infinite cylinder which may have an oscillating boundary and/or a periodic array of holes.
\end{remark}

{\bf Proof of Proposition \ref{phodge}.} We first consider the case where $d=2$. We define
$$
\hat\Omega := \Omega/(L_1\mathbb{Z} \times L_2\mathbb{Z})~~\text{ and }~~
T := \R^2/(L_1\mathbb{Z} \times L_2\mathbb{Z}).
$$

If $x \in \R^2$, we denote by $\hat x$ its class of equivalence in $T$, and if $\phi:\R^2 \to \R$ is $L$- periodic, we denote $\hat\phi $ the function $T \to \R^2$ verifying $\phi(x)=\hat\phi(\hat x)$.
\vskip 0.3cm

Finding $\phi \in C^{2,\delta}(\overline{\Omega}),$ which is $L$-periodic with respect to $x$ and verifying (\ref{hodge}), is then equivalent to finding $\hat \phi \in C^{2,\delta}(\overline{\hat \Omega})$ verifying (\ref{hodge}). We consider the map $\tilde{ q}$ defined as follows:
\begin{eqnarray*}
\tilde{q} : T & \longrightarrow & \R^2, \\
 x\in \overline{\hat \Omega} & \longmapsto & q(x), \\
\hat x \notin \overline{\hat{\Omega}} & \longmapsto & 0.
\end{eqnarray*}
We claim that $\tilde{q}$ is a divergence free vector field on $T$ in the sense of distributions. Indeed, if $\psi \in C^\infty(T),$ we then have
\begin{eqnarray*}
< div(\tilde{q}), \psi> & := & -<\tilde{q},\nabla \psi>~=~-\int_T \tilde q \cdot \nabla \psi \vspace{3 pt}\\
& = & \ds{-\int_{\hat{\Omega}} q \cdot \nabla \psi} =\ds{-\int_{\partial \hat{\Omega}} \psi\, q \cdot \nu + \int_{\hat{\Omega}} \psi \nabla \cdot q } \vspace{3 pt} \\
& = & 0 + 0 = 0,
\end{eqnarray*}
because of the conditions (\ref{q}). Moreover, we clearly have $\displaystyle \int_T \tilde{q}=0$.
Now, we denote by $R$ the matrix of a direct rotation of angle $\pi/2$.

The next step is to solve the following equation in $H^1(T)$:
\begin{equation}\label{laplace}
-\Delta \tilde{\phi} = \nabla \cdot (R\tilde{q}).
\end{equation}

A function $\tilde{\phi} \in H^1(T)$ is a weak solution of (\ref{laplace}) if we have, for all $\psi \in H^1(T),$
\begin{equation}\label{wlaplace}
\int_T \nabla \tilde{\phi} \cdot \nabla \psi = - \int_T R \tilde{q} \cdot \nabla \psi.
\end{equation}
We set $\displaystyle E:=\{ \psi \in H^1(T) \text{ such that } \int_T \psi = 0\}$, so that, thanks to Poincar\'e's inequality,
$$
<u,v>_E := \int_T \nabla u \cdot \nabla v
$$
is an inner product on $E$. Moreover, $\displaystyle \psi \in E \mapsto \int_T R \tilde{q} \cdot \nabla \psi$ is a continuous linear form on $E$, so by the Lax-Milgram theorem, there exists a unique $\tilde{\phi} \in E$ solution of (\ref{wlaplace}).
The condition $\int_T \psi=0$ is not restrictive because only the gradients of functions belonging to $E$ appear in the weak formulation (\ref{wlaplace}). We then have $\tilde{\phi} \in H^1(T)$ such that in the sense of distributions
\begin{eqnarray*}
\begin{array}{lll}
\nabla \cdot R (\tilde{q} - \nabla^\perp \tilde{\phi})& = & 0 \text { in }T\quad\hbox{and}\\
\nabla \cdot (\tilde{q} - \nabla^\perp \tilde{\phi}) & = & 0 \text { in }T~~\text{ since $\nabla\cdot\tilde{q}=0$ in $\mathcal{D}'(T)$ and $div(\nabla^\perp\cdot)=0$}.
\end{array}
\end{eqnarray*}
This implies that $\tilde q - \nabla^\perp \tilde{\phi}$ is a harmonic distribution on $T$. Using Weyl's theorem (see \cite{J Just}), we conclude that $\tilde q - \nabla^\perp \tilde{\phi}$ is a harmonic function on the torus $T$ and therefore is constant. Indeed, if $h$ is a harmonic scalar function on $T$, then by multiplying $h$ by $\Delta h$ and integrating by parts we get
$$
\int_T|\nabla h|^2 = 0,
$$
which leads to $h$ is constant on $T$.

Finally, since $\displaystyle \int_T (\tilde{q} - \nabla^\perp \tilde{\phi}) = 0$ and $\tilde{q} - \nabla^\perp \tilde{\phi}$ is constant, we conclude that, in the sense of distributions:
\begin{equation}\label{dhodge}
\tilde{q} = \nabla^\perp \tilde{\phi}.
\end{equation}
We set $\hat \phi := \tilde{\phi}|_{\hat{\Omega}}$, which solves (\ref{hodge}) in $\hat{\Omega}$. The corresponding $L$-periodic function $\phi \in H^1_{loc}(\Omega)$ solves then (\ref{hodge}) in $\Omega$.   The $C^{2,\delta}$ regularity of $\phi$ in $\overline{\Omega}$ is a consequence of the Schauder estimates for the Laplace equation.

The fact that $\phi$ is constant on every connected component of $\partial \Omega$ is a straightforward consequence of the identity
$$
\nabla^\perp \phi \cdot \nu = q \cdot \nu = 0 \text{ on } \partial \Omega.
$$

For the case $d=1$, we symmetrize the set $\Omega$ (resp. the cell $C$) and the field $q$ with respect to the line $y=R$ and we call the resulting set by $\Omega_s$ (resp. $C_s$) and the resulting
vector field by $q_s.$ For the sake of completeness, we mention that $q_s(x,y)$ is given by
$$q_s(x,y)=\left\{ \begin{array}{ll}
q(x,y)\hbox{ for } (x,y)\in\Omega,\\
q_1(x,2R-y)e_1\mathbf{-}q_2(x,2R-y)e_2\hbox{ for }(x,y)\in \Omega_s\setminus \Omega.
\end{array}\right.$$
One can easily notice that $\int_{C_s}q_s=0.$
 We then generate (in the direction $e_2$) a \textbf{periodic set} $\Omega_1$ from the set $\Omega_s$  in order to reduce this case to the case $d=2$. For that purpose, since we already have
$
\Omega_s \subset \R \times [-R,3R]
$
(take $d=1$ and $N=2$ in (\ref{comega})),
we define $\Omega_1$ in the following way
\begin{equation}\label{omega}
\Omega_1 := \bigcup_{i \in \mathbb{Z}} \ds{\left\{\Omega_s + i(4R+2)e_2\right\}}.
\end{equation}
Thus, $\Omega_1$ is periodic in the direction of $e_1$ and $e_2$ and is the disjoint reunion of translations of $\Omega_s$. We set
$$
\hat{\Omega}:=\Omega_1 / (L_1\mathbb{Z} \times (4R+2)\mathbb{Z}) \ \text{ and } \ T:=\R^2/ (L_1\mathbb{Z} \times (4R+2)\mathbb{Z}),
$$
and the procedure used for the case $d=2$ still works in this case and gives $\phi \in C^{2,\delta}(\overline{\Omega})$, $L$-periodic, solving
$$
q = \nabla^\perp \phi.
$$\hfill$\Box$
\begin{corollary}\label{firstint}
Let
\begin{equation}\label{j}
\mathcal{J}:=\left\{\eta \circ \phi,\text{ such that } \eta : \R \to \R \ \text{ is Lipschitz} \right\},
\end{equation}
where $\phi$, such that $q = \nabla^\perp \phi$, is given by Proposition \ref{phodge}. Then,
$$
\mathcal{J} \subset \mathcal{I}\cup\{0\}.
$$
\end{corollary}
{\bf Proof.} We first mention that if $q \not\equiv 0$ is the advection vector field and $\phi$ is the function given by Proposition \ref{phodge}, then $\phi \in \mathcal{I}$. Indeed, $\phi \in H^1_{loc}(\Omega)$ and we have
$$
\forall\, z \in \Omega,~~q(z) \cdot \nabla \phi(z) = \nabla^\perp \phi(z) \cdot \nabla \phi(z) = 0.
$$
Now, using remark \ref{about I}, we conclude that $\eta \circ \phi$ is a first integral whenever $\eta :\R \to \R$ is a Lipschitz function.\hfill$\Box$

\begin{remark}
For any $w \in \mathcal{J}$, we have $\int_C (q \cdot \tilde{e})w^2 = 0.$ Indeed, $w = \eta \circ \phi$ and $q = \nabla^\perp \phi$ which gives
\begin{eqnarray*}
\int_C (q \cdot \tilde{e})w^2 & = & \tilde{e} \cdot \int_C \left(\nabla^\perp \phi\right)\, \eta^2(\phi) \\
& = & \tilde{e} \cdot R \int_C \nabla \left(F\circ\phi\right),
\end{eqnarray*}
where $R$ the matrix of a direct rotation of angle $\pi/2$ and $F' = \eta^2.$ Let $\tilde{\phi} \in H^1(T)$ defined by (\ref{dhodge}), then $\tilde{\phi} = \phi$ on $\hat{\Omega}$ and is constant on every connected component of $T \backslash \hat{\Omega}$, and so is $\ds{F\circ\tilde{\phi}}.$ We then have
$$
\ds\int_{\ds{T\backslash \ds\hat{\Omega}}} \nabla \left(F\circ\tilde{\phi}\right)= 0.
$$
Hence,
\begin{eqnarray*}
\int_C (q \cdot \tilde{e})w^2 & = & \tilde{e} \cdot R \int_T \nabla \left(F\circ\tilde{\phi}\right)=0,
\end{eqnarray*}
because $T$ has no boundary.
Thus,
\begin{equation}\label{jzero}
\forall \,w \in \mathcal{J},~~\int_C (q \cdot \tilde{e})w^2 = 0.
\end{equation}
We recall that the family of first integrals $\mathcal{I}$ always contains the set $\mathcal{J}.$ However, this does not, in general, provide enough information about the following quantities
$$
\sup_{w \in \mathcal{I}} \int_C (q \cdot \tilde{e})w^2 \text{ or } \max_{w \in \mathcal{I}_1} \int_C(q \cdot \tilde{e})w^2,
$$
which appear in the asymptotics of the minimal speed within a large drift.
Lemmas \ref{etaphi} and \ref{stephane} are devoted to prove Theorem \ref{case N=2} and treat this situation.
\end{remark}
\begin{definition}
Throughout the rest of this section, we denote
\begin{equation} \label{omega}
\left\{
\begin{array}{lllllllll}
T & := & \R^2/(L_1 \Z \times L_2 \Z)& \hbox{and} & \hat \Omega & := & \Omega/(L_1 \Z \times L_2 \Z) & \text{if} & d=2,\vspace{4 pt} \\
T & := & \R^2/\left(L_1\Z\times\{0\}\right) & \hbox{and}  & \hat \Omega & := & \Omega/\left(L_1\Z\times \{0\}\right) & \text{if} & d=1.
\end{array}
\right.
\end{equation}
Moreover, if $x \in \Omega$ (resp. $\R^2$), $\hat x$ denotes its class of equivalence in $\hat \Omega$ (resp. $T$), and if $u:\Omega \to \R$ is $L$-periodic, $\hat u$ denotes the function $\hat \Omega \to \R$ verifying $\hat u (\hat x)=u(x)$ for almost every $x \in \Omega$.\\

We also define the canonical projection on $T$ by
\begin{eqnarray}\label{proj}
\Pi : \R^2 & \longrightarrow & T \\
x & \longmapsto & \hat x. \nonumber
\end{eqnarray}
\end{definition}

We need the following preliminary lemma in order to prove the main theorem of this section:
\begin{lemma}\label{etaphi}
Let $\hat{\Omega}$ be the set defined in (\ref{omega}), $\hat U$ be an open subset of $\hat{\Omega}$, and $\hat \phi$ given by (\ref{hodge}). We suppose the following:

(i) $\hat q(\hat x) \neq 0$ for all $\hat x \in \hat U,$

(ii) the level sets of $\hat \phi$ in $\hat U$ are all connected.\\
Then, for every $w \in \mathcal{I}$, there exists a continuous function $\eta: \hat \phi( \hat U) \to \R$ such that
\begin{equation}\label{eqetaphi}
\hat w = \eta \circ \hat \phi \text{ on }  \hat U.
\end{equation}
\end{lemma}
{\bf Proof.} For every $\lambda \in \hat \phi(\hat U)$, we denote by $\Gamma_\lambda$ the level set $$\Gamma_\lambda:=\{ \hat x \in \hat U \text{ such that }\hat \phi(x)= \lambda\}.$$ It follows from (i) that $\nabla \hat \phi$ does not vanish on $\hat U$; and hence, $\hat w$ has to be constant on the connected components of the level sets of $\hat \phi$ because
$$
\forall \hat x \in \hat U, \ \nabla^\perp \hat \phi(\hat x) \cdot \nabla \hat w(\hat x)=q(x)\cdot\nabla w(x)=0.
$$
By (ii), the level sets of $\hat \phi$ in $\hat U$ are connected, so $\hat w$ is constant on every level set of $\hat\phi$. If $\hat x \in \Gamma_\lambda$, we have $\hat\phi(\hat x)=\lambda$, and $\hat w$ is constant on $\Gamma_\lambda$, so depends only on $\lambda$. Then, we can define $\eta$ by
$$
\eta(\lambda) = \hat w(\Gamma_\lambda).
$$
To prove the continuity of $\eta,$ we suppose, to the contrary, that there exists $\lambda_0 \in \hat\phi(\hat U)$ such that $\eta$ is not continuous at $\lambda_0$.  The set $\Gamma_{\lambda_0}$ is a curve because $\hat U$ is open and $\nabla \hat \phi$ does not vanish on $\hat U$ by (i). The function $\hat w$ then has a ``jump'' along the level set $\Gamma_{\lambda_0}$, which is impossible because $\hat w \in H^1(\hat U)$, and so has a trace on $\Gamma_{\lambda_0}$.\hfill $\Box$
\vskip0.5cm

Recalling Definition \ref{trajectory def} of the trajectories of an advection field $q,$ we mention that:
\begin{remark}
It is obvious that $\phi$ is constant on every trajectory of $q$. Moreover,
the trajectories of $q$ make a partition of $\Omega \backslash \{ x \in \Omega$ such that $q(x)=0\}$.\\
If $T(x)$ is the trajectory of $q$ in $\Omega$, and $T(\hat x)$ is the trajectory of $\hat q$ in $\hat \Omega$, then we have
$$
T(\hat x) =\Pi(T(x)).
$$
In other words, the trajectory of the projection is the projection of the trajectory.
\end{remark}
\begin{definition}\label{def} We define here the set of ``regular trajectories'' in $\hat{\Omega}$. Let
\begin{equation}\label{U hat}
\hat U:=\left\{\hat x \in \hat \Omega\text{ such that }T(\hat x)\text{ is well defined and closed in } \overline{\hat \Omega}\right\}.
\end{equation}
We denote by $\ds{\hat U_i}$ the connected components of $\ds{\hat U}$.
\end{definition}
The set $\hat{U}$ is exactly the union of the trajectories which are simple closed curves in $\hat \Omega$. This is proved in the following proposition:
\begin{proposition}\label{trajec}
Let $\hat x \in \hat U$, then $T(\hat x)$ is a $C^1$ simple closed curve in $\overline{\hat \Omega}$.
\end{proposition}
{\bf Proof.} Let $\hat x \in \hat U$. By definition of $\hat U$, $T(\hat x)$ is closed. Moreover, since $T(\hat x)$ is a subset of $\overline{\hat \Omega}$ which is a compact set, $T(\hat x)$ is then compact. $\hat x \mapsto |\hat q (\hat x)|$ attains then its minimum on $T(\hat x)$ at some point $\hat x_0 \in T(\hat x)$. Since $\hat x_0 \in T(\hat x)$ we have $|\hat q (\hat x_0)|= \eta > 0$. We then get
\begin{equation}\label{phiq}
\min_{\hat y \in T(\hat x)} |\hat q(\hat x)| \geq \eta >0.
\end{equation}
Besides, we know that $\hat \phi$ is constant on $T(\hat x)$. Let $\alpha := \hat \phi(T(\hat x))$, and $A_\alpha := \{\hat x \in \hat \Omega$ such that $ \phi(x) > \alpha \}$. Since $\nabla \hat \phi$ does not vanish on $T(\hat x)$ because of (\ref{phiq}), we get $T(\hat x) \subset \partial A_\alpha$. Using the Stokes formula on $A_\alpha$ with the vector field $\nabla \hat \phi$ gives
$$
\int_{A_\alpha} -\Delta \hat \phi = \int_{\partial A_\alpha} | \nabla \hat\phi|.
$$
Thus,
$$
\eta \mathcal{L}^1(T(\hat x)) \leq \int_{T(\hat x)} | \nabla \hat \phi| \leq \int_{\hat \Omega} |\Delta \hat \phi| < +\infty,
$$
where $\mathcal{L}^1$ denote the 1-dimensional Lebesgue measure.
The trajectory $T(\hat x)$ is then a $C^1$ curve with finite length. It has no self intersection point because such a point would be a critical point of $\hat \phi$. $T(\hat x)$ has no boundary because, if $\hat y \in \partial T(\hat x)$, since $T(\hat x)$ is closed, $\hat y \in T(\hat x)$, and so $\hat q (\hat y) \neq 0$ and we could extend $T(\hat x)$ at the point $\hat y$.

Consequently, if $\hat x \in \hat U$, $T(\hat x)$ is a $C^1$ simple closed curve with finite length.\hfill$\Box$

\begin{lemma}\label{stephane}
Let $\ds{\hat U_i}$ as in the previous definition. Then,

(i) all the level sets of $\hat \phi$ in $\hat U_i$ are connected,

(ii) all the level sets of $\hat \phi$ in $\hat U_i$ are homeomorphic,

(iii) $\partial \hat U_i$ has exactly two connected components $\hat \gamma_1$ and $\hat \gamma_2$ such that
$$
\hat \phi(\hat \gamma_1) = \sup_{\hat x \in \hat U_i} \hat \phi(\hat x) \text{ and } \hat \phi(\hat \gamma_2) = \inf_{\hat x \in \hat U_i} \hat \phi(\hat x).
$$
\end{lemma}

{\bf Proof.} Let $\Gamma_\lambda :=\{\hat x \in \hat U_i$ such that $\hat \phi(\hat x)=\lambda\}$ be a non empty level set of $\hat \phi$ in $\hat U_i$. Let $\Gamma_\lambda^1$ be one of its connected components. $\Gamma_\lambda^1$ is a $C^1$ curve because $\hat \phi \in C^{2,\delta}$ and $\nabla \hat \phi$ does not vanish on $\Gamma_\lambda^1$ by definition of $\hat U$. Let $\hat x \in \Gamma_\lambda^1$, we consider the following ODE:
\begin{equation}\label{cauchy}
\left\{
\begin{array}{lll}
y'(t) & = & \ds{\frac{\nabla \hat \phi (y(t))}{|\nabla \hat \phi(y(t))|^2}},\vspace{4 pt}\\
y(0) & = & \hat x.
\end{array}
\right.
\end{equation}
By classical ODE theory, there exists a maximal interval $\ds{(t^1_{\hat x},t^2_{\hat x})}$ with $t^1_{\hat x}<0<t^2_{\hat x}$ on which there exists a $C^1$ solution $y_{\hat x}$ of (\ref{cauchy}). Moreover, either $y(t^1_{\hat x})$ (resp. $y(t^2_{\hat x})$) is a critical point of $\hat \phi$ or belongs to $\partial \hat{\Omega}$, otherwise we could extend the maximal solution $y_{\hat x}$ to a larger interval, which is impossible. We set
$$
\ds{t^1:=\max_{\hat x \in \Gamma_\lambda^1} t^1_{\hat x}} ~~\text{ and }~~\ds{ t^2:=\min_{\hat x \in \Gamma_\lambda^1} t^2_{\hat x}}.
$$
Since $t^1$ (resp. $t^2$) is equal to $t^1_{\hat x}$ (resp. $t^2_{\hat x}$) for some point $\hat x \in \Gamma_\lambda^1,$ we then have
$$
t_1 < 0 < t_2.
$$
We define the function $g$ on $(t_1,t_2) \times \Gamma_\lambda^1$ by
$$
g_t(\hat x) = y_{\hat x}(t).
$$
We claim that $\phi(g_t(\hat x)) = \lambda+t$ for every $(t,\hat x) \in (t_1,t_2)\times \Gamma_\lambda^1$. Indeed, we have
\begin{eqnarray*}
\frac{d}{dt}\hat \phi(g_t(\hat x)) & = & \ds{\frac{d}{dt}\hat \phi(y_{\hat x}(t))}\vspace{3pt} \\
& = & \ds{y'_{\hat x}(t) \cdot \nabla \hat \phi(y_{\hat x}(t)) }\vspace{3pt}\\
& = & \ds{\frac{\nabla \hat \phi (y_{\hat x}(t))}{|\nabla \hat \phi(y_{\hat x}(t))|^2} \cdot \nabla \hat \phi(y_{\hat x}(t))}\vspace{3pt} \\
& = & 1,
\end{eqnarray*}
which leads to $\hat \phi(g_t(\hat x)) - \hat \phi(\hat x) = t = \hat \phi(g_t(\hat x)) - \lambda$. Thus, on $g_t(\Gamma_\lambda^1)$, $\hat \phi$ is  equal to the constant $\lambda+t$. Moreover, $g_t$ is continuous with respect to $\hat x$ due to the continuity of the solution of an ODE with respect to initial data. We set now
$$
\hat V_i := \bigcup_{t_1<t<t_2} g_t(\Gamma_\lambda^1).
$$
We need the following two claims,  whose proofs are postponed at the end, in order to prove that $\hat V_i=\hat U_i:$ \vskip 0.35cm

{\bf Claim 1:} For every $\varepsilon >0$, there exists $r_\varepsilon$ such that for every $\hat x \in \hat V_i,$
$$
\hat x \in \bigcup_{t_1+\varepsilon \leq t \leq t_2-\varepsilon}g_t(\Gamma_\lambda^1) \Rightarrow B(\hat x,r_\varepsilon) \subset \hat V_i,
$$
and as a consequence $\hat V_i$ is an open subset of $T$.\vskip0.35cm

{\bf Claim 2:} $\partial \hat V_i$ has exactly two connected components $C_1$ and $C_2$ such that
$$
\hat \phi|_{C_1} \equiv \lambda+t_1 \text{ and } \hat \phi|_{C_2} \equiv \lambda+t_2,
$$
and either $C_1$ (resp. $C_2$) is a connected component of $\partial \hat{\Omega}$ or contains a critical point of $\hat \phi$.\\

 By definition, $\hat V_i$ is the union of connected components of level sets of $\hat \phi$ on which $\nabla \hat \phi$ does not vanish. Hence, $\hat V_i \subset \hat U.$ Moreover, $\hat V_i$ is connected because, by construction, any point $\hat x \in \hat V_i$ can be connected to $\Gamma_\lambda^1$, which is a connected set. Finally, $\hat U_i \cap \hat V_i \neq \emptyset$ because it contains $\Gamma_\lambda^1$. We then can affirm that $\hat V_i \subset \hat U_i$.

Suppose that this inclusion is strict, then we can find $\hat x_0 \in \hat U_i \backslash \hat V_i$. Let $\hat x_1 \in \hat V_i$, and $\gamma:[0,1] \to \hat U_i$ a continuous path connecting $\hat x_0$ and $\hat x_1$. By continuity, it crosses $\partial \hat V_i.$ However, $\partial \hat V_i \cap \hat U_i=\emptyset$ because, by claim 2, the connected components of $\partial \hat V_i$ are either connected components of $\partial \hat{\Omega}$, which do not intersect $\hat U_i$, or contain a critical point of $\hat \phi$, and are thus removed from $\hat U$ by construction.

Properties (i) and (ii) are straightforward, because a level set of $\hat U_i$ is a level set of $\hat V_i$ and can be written in the form $g_t(\Gamma_\lambda^1)$, and $g_t$ is a homeomorphism for any $t_1 < t < t_2$. The level sets are then all homeomorphic to $\Gamma_\lambda^1$, which is connected.

Property (iii) has already been proved for $\hat V_i$, and $\hat V_i=\hat U_i$. Eventually, the proof of the lemma is complete.\hfill$\Box$
\vskip 0.35cm

{\bf \underline{Proof of claim 1.}} We prove this claim by contradiction. First, we have
$$
\bigcup_{t_1+\varepsilon \leq t \leq t_2-\varepsilon}g_t(\Gamma_\lambda^1) = g([t_1+\varepsilon,t_2-\varepsilon],\Gamma_\lambda^1),
$$
where we denote $g(t,\hat x)=g_t(\hat x)$. Moreover, $g$ is a continuous function, because of the continuity of the solution of an ODE with respect to the initial conditions. Therefore, $\ds{\bigcup_{t_1+\varepsilon \leq t \leq t_2-\varepsilon}g_t(\Gamma_\lambda^1)}$ is compact, as the image of a compact by a continuous mapping.

Suppose now that claim 1 is not true. Then, there exist $\varepsilon>0$ and $\hat x$ such that $\hat x \in \partial \hat V_i$ and $\hat \phi (\hat x)=: \lambda + \alpha \in [\inf_{\hat V_i} \hat \phi + \varepsilon,\sup_{\hat V_i} \hat \phi - \varepsilon ]$.

Let $\psi$ defined on $(-\beta,\beta)$ be an arc length local parametrization of $g_\alpha(\Gamma_\lambda^1)$ such that $\psi(0) = \hat x$. We define the following mapping for $0<\xi<\varepsilon$:
\begin{eqnarray*}
G :(-\beta,\beta) \times (-\xi,\xi) & \longrightarrow & \hat V_i \\
(s,t) & \longmapsto & g_{t}(\psi(s)).
\end{eqnarray*}
We have $G(0,0)=\hat x$ and
$$
DG(0,0)(s,t) = s \psi'(0) + t \frac{\nabla \hat \phi (\hat x)}{|\nabla \hat \phi (\hat x)|^2}.
$$
The linear mapping $DG(0,0)$ is then an automorphism of $\R^2$, because $$\ds{\left\{\psi'(0),\frac{\nabla \hat \phi (\hat x)}{|\nabla \hat \phi (\hat x)|^2}\right\}}$$ is an orthogonal basis of $\R^2$. The application of the inverse mapping theorem then gives two open sets $W_1 \subset (-\beta,\beta) \times (-\xi,\xi)$, with $(0,0) \in W_1,$ and $W_2 \subset \hat V_i,$ with $\hat x \in W_2$, such that $G$ is a local diffeomorphism form $W_1$ to $W_2$. This prevents $\hat x$ from belonging to $\partial \hat V_i$ and then gives a contradiction.

The fact that $\hat V_i$ is an open subset of $T$ is a straightforward consequence, because if $\hat x \in \hat V_i$, then for $\varepsilon$ small enough $\hat x \in \bigcup_{t_1+\varepsilon \leq t \leq t_2-\varepsilon}g_t(\Gamma_\lambda^1)$ and so we can find a neighborhood of $\hat x$ in $\hat V_i$.
\vskip 0.35cm

{\bf \underline{Proof of claim 2.}} Let $\hat x \in \partial \hat V_i$, then $\hat \phi (\hat x)=\lambda + t_1$ or $\lambda + t_2$. indeed, by continuity of $\hat \phi$ we have $\hat \phi (\hat x) \in [t_1,t_2]$. Suppose by contradiction that $\hat \phi (\hat x) = \alpha \in (t_1,t_2)$, then for $\varepsilon$ sufficiently small we have $\alpha \in [t_1+2\varepsilon,t_2-2\varepsilon]$. Let $\{\hat x_p \}$ be a sequence in $\hat V_i$ converging to $\hat x$. We have then $\hat \phi(\hat x_p) \to \alpha$ as $p \to \infty$. Hence for $p$ large enough we have $\hat \phi(\hat x_p) \in [t_1+\varepsilon,t_2-\varepsilon]$, so by claim 1, there exists $r>0$ such that for $p$ large enough $dist(\hat x_p,\partial \hat V_i)\geq r>0$, leading to $dist(\hat x,\partial \hat V_i)\geq r>0$ which contradicts the fact that $\hat x \in \partial \hat V_i$.\\

We now prove that $\partial \hat V_i$ has exactly two connected components. Let $C_1:=\partial \hat V_i \cap \hat \phi^{-1}(\lambda+t_1)$ and  $C_2:=\partial \hat V_i \cap \hat \phi^{-1}(\lambda+t_2)$. We have, using the previous remark,
$$
C_1 = \bigcap_{p \in \mathbb{N}} \overline{\hat V_i \cap \hat \phi^{-1}((\lambda + t_1,\lambda + t_1+1/p))}.
$$
Since $\hat V_i \cap \hat \phi^{-1}((\lambda + t_1,\lambda + t_1+1/p))$ is a nonempty bounded open connected subset of $T$, its closure is a nonempty compact connected subset of $T$. Therefore $C_1$ is the decreasing intersection of nonempty connected compact subsets of $T$, and is then a connected nonempty compact subset of $T$. Similarly, $C_2$ is connected.\\

Finally, we prove that $C_1$ (resp. $C_2$) is either a connected component of $\partial \hat \Omega$ or contains a critical point of $\hat \phi$. We suppose then that $C_1$ does not contain any critical point of $\hat \phi$, it must then contain a point $\hat x_0$ of $\partial \hat \Omega$, otherwise for any $\hat x \in \Gamma_\lambda^1$ the solution of (\ref{cauchy}) could be extended at the point $t_1$, and this would contradict the definition of $t_1$. We denote by $D_1$ the connected component of $\partial \hat \Omega$ containing $\hat x_0$. We are left to prove that $D_1=C_1$. First, we know that $C_1$ is a regular simple closed curve, because it is a connected component of a level set of $\hat \phi$ on which $\hat \phi $ does not vanish. We have then $C_1=T(\hat x_0)$. Moreover the trajectories of $\hat q$ intersecting the boundary of $\hat \Omega$ follow the boundary of $\hat \Omega$ since $q\cdot\nu=0$ on $\partial \Omega$, so $T(\hat x_0) \subset D_1$. We conclude that $C_1 \subset D_1$, and since $C_1$ and $D_1$ are both connected simple closed curves we get $C_1=D_1$. Similarly, we get that $C_2$ contains a critical point of $\hat \phi$ or is a connected component of $\partial \hat \Omega$.\hfill$\Box$

\vskip0.5cm

{\bf Proof of Theorem \ref{case N=2}.} Using (\ref{hodge}), we have for any $w \in \mathcal{I}$
$$
\int_C q w^2 =R \int_C (\nabla \phi) w^2 = R \int_{\hat{\Omega}} (\nabla \hat \phi)\hat w^2.
$$
Let $W:=\{ \hat x \in \hat{\Omega}$ such that $\hat \phi(\hat x)$ is a critical value of $\hat \phi\}$. Using the co-area formula (\cite{EG}, \cite{Federer}) we get
$$
\left|\int_W \hat w^2 \nabla \hat \phi\right| \leq \int_W \hat w^2 |\nabla \hat \phi| = \int_{\hat \phi(W)} \left(\int_{\hat \phi^{-1}(t)}\hat w^2(x)\right)dt .
$$
Moreover, from Sard's theorem (see \cite{sard} for eg.), since $\hat \phi$ is $C^2$, $\mathcal{L}^1(\hat \phi(W))=0$, where $\mathcal{L}^1$ denotes the Lebesgue measure on $\R.$ It follows that
$$
\int_W \hat w^2 \nabla \hat \phi = 0.
$$
Since $\hat{\Omega} \backslash W \subset \hat U \subset \hat{\Omega},$ we get
\begin{equation}\label{split}
\int_C q w^2 = R \int_{\hat U} (\nabla \hat \phi) \hat w^2 = R \sum_i \int_{\hat U_i}(\nabla \hat \phi) \hat w^2.
\end{equation}
We now use Lemma \ref{etaphi} to get $\eta_i$ continuous such that
$$
\int_{\hat U_i}(\nabla \hat \phi) \hat w^2 = \int_{\hat U_i}(\nabla \hat \phi) \eta_i^2(\hat \phi).
$$
We define the function $F_i$ by $F_i'=\eta_i^2$ and $F_i(0)=0$, and we obtain
$$
\int_{\hat U_i}(\nabla \hat \phi) \hat w^2 = \int_{\hat U_i}\nabla F_i(\hat \phi).
$$
If we define
$$
\hat U_i^\varepsilon :=\{ \hat x \in \hat U_i \text{ such that }\inf_{\hat U_i}\hat \phi + \varepsilon < \hat \phi(x) < \sup_{\hat U_i} \hat \phi - \varepsilon\},
$$
then it follows from dominated convergence theorem that
\begin{equation}\label{domin}
\int_{\hat U_i^\varepsilon} (\nabla \hat \phi) \hat w^2 \xrightarrow[\varepsilon \to 0]{} \int_{\hat U_i}(\nabla \hat \phi) \hat w^2.
\end{equation}

We now prove (i) $\Rightarrow$ (ii) by contraposition.  We suppose that there exist no periodic unbounded trajectories of $q$. In $\hat U_i$, the trajectories of $q$ are exactly the level sets of $\hat \phi$. We consider the following set
$$
U_i^\varepsilon:= \Pi^{-1}(\hat U_i^\varepsilon).
$$
Let $x_0 \in U_i^\varepsilon$ and let $U_{i,0}^\varepsilon$ be the connected component of $U_i^\varepsilon$ containing $x_0$. We claim that $\Pi$ is a  surjection from $U_{i,0}^\varepsilon$ to $\hat U_i^\varepsilon$. For that purpose we prove that $\Pi(U_{i,0}^\varepsilon)$ is open and closed in $\hat U_i^\varepsilon$.

Let $\hat x \in \Pi(U_{i,0}^\varepsilon)$, and $r>0$ sufficiently small to have $B(\hat x,r) \subset \hat U_i^\varepsilon$. Let $x \in U_{i,0}^\varepsilon$ such that $\Pi(x)=\hat x$, then $\Pi(B(x,r)) = B(\hat x,r)$ which leads to $B(\hat x,r) \subset \Pi(U_{i,0}^\varepsilon)$, so $\Pi(U_{i,0}^\varepsilon)$ is open in $\hat U_i^\varepsilon$.

On the other hand, let $\{\hat x_n \}$ be a sequence of $\Pi(U_{i,0}^\varepsilon)$ converging to $\hat x \in \hat U_i^\varepsilon$. Since $\hat U_i^\varepsilon$ is open then, for $r>0$ sufficiently small, $B(\hat x,r) \subset \hat U_i^\varepsilon$. Hence, for $n$ large enough, we have $\hat x_n \in B(\hat x,r)$, so there exists $N \in \mathbb{N}$ and $s>0$ such that $B(\hat x_N,s)\subset B(\hat x,r)$ and $\hat x \in B(\hat x_N,s)$. Finally, let $x_N \in U_{i,0}^\varepsilon$ such that $\Pi(x_N) = \hat x_N$, then there exists $x \in B(x_N,s)$ such that $\Pi(x) = \hat x$, which leads to $\Pi(U_{i,0}^\varepsilon)$ is closed in $\hat U_i^\varepsilon$.

Now, by definition, $\hat U_i$ and $\hat U_i^\varepsilon$ only contain ``regular'' trajectories of $q$, and so does $U_{i,0}^\varepsilon$. By assumption, there exist no periodic unbounded trajectories of $q$, so all the trajectories of $q$ in $U_{i,0}^\varepsilon$ are bounded. Moreover, $\partial U_{i,0}^\varepsilon$ is the disjoint reunion of two bounded, regular, trajectories of $q$. All the trajectories of $q$ in $U_{i,0}^\varepsilon$ are level sets of $\phi$, and, by a compactness argument we get that $U_{i,0}^\varepsilon$ is bounded in $\Omega$.

From this boundedness, we obtain that $\Pi:U_{i,0}^\varepsilon \to \hat U_i^\varepsilon$ is a bijection. Indeed, if it is not injective, since $U_{i,0}^\varepsilon$ is connected we would have a path connecting two different points $x_1$ and $x_2$ of $U_{i,0}^\varepsilon$, such that $\Pi(x_1) = \Pi(x_2)$, and by periodicity, $U_{i,0}^\varepsilon$ could not be bounded.

We conclude that $\Pi:U_{i,0}^\varepsilon \to \hat U_i^\varepsilon$ is a measure preserving bijection, by definition of the measure on $T$. We get then
\begin{eqnarray*}
\int_{\hat U_i^\varepsilon} (\nabla \hat \phi) \hat w^2 & = & \int_{U_{i,0}^\varepsilon} (\nabla \phi) w^2\\
& = &  \int_{U_{i,0}^\varepsilon} \nabla F_i(\phi)\\
& = & \int_{\partial U_{i,0}^\varepsilon} F_i(\phi) \mathbf{n},
\end{eqnarray*}
 where $\mathbf{n}$ is the unit outward normal vector field to $\partial U_{i,0}^\varepsilon$. Finally $\partial U_{i,0}^\varepsilon$ is the union of two level sets $C_1$ and $C_2$ of $\phi$ in $\Omega$, which are both simple closed curves, so we can write
$$
\int_{U_{i,0}^\varepsilon} (\nabla \phi) w^2 = F_i(\phi(C_1))\int_{C_1}\mathbf{n} + F_i(\phi(C_2))\int_{C_2}\mathbf{n},
$$
with
$$
\int_{C_1}\mathbf{n} = \int_{C_2}\mathbf{n} = 0,
$$
because the integral of the unit normal on a $C^1$ closed curve in $\R^2$ is zero.
Therefore, using (\ref{split}) and (\ref{domin}), we get
$$
\int_C q w^2=0,
$$
for all $w \in \mathcal{I}$. \vskip 0.5cm

We now prove (ii) $\Rightarrow$ (i). Let $x \in \Omega$, and $\mathbf{a}\neq 0$, such that $T(x)$ is $\mathbf{a}-$periodic and unbounded such that $|\mathbf{a}|$ is minimal. Let $\hat U_i$ be the connected component of $\hat U$ containing $\hat x := \Pi(x)$. We define as previously $\hat U_i^\varepsilon$ for $\varepsilon$ sufficiently small in order to have $\hat U_i^\varepsilon \neq \emptyset$. Let $x_0 \in \Pi^{-1}(\hat U_i^\varepsilon)$, we define $U_{i,0}^\varepsilon$ to be the connected component of $\Pi^{-1}(\hat U_i^\varepsilon)$ containing $x_0$. Let $e'_1 := \mathbf{a}/|\mathbf{a}|$, and $e'_2$ such that $(e'_1,e'_2)$ is an orthonormal frame of $\R^2$.
We set
$$
U_i^\varepsilon := \{ x=x'_1e'_1+x'_2e'_2 \in U_{i,0}^\varepsilon \text{ such that } 0\leq x'_1 <|\mathbf{a}| \}.
$$
Using similar arguments to (i) $\Rightarrow$ (ii), we get that $\Pi$ is a measure preserving bijection from $U_i^\varepsilon$ to $\hat U_i^\varepsilon$, and we have
\begin{eqnarray*}
\int_{\hat U_i^\varepsilon} (\nabla \hat \phi) \hat w^2 & = & \int_{U_i^\varepsilon} (\nabla \phi) w^2\\
& = &  \int_{U_i^\varepsilon} \nabla F_i(\phi)\\
& = & \int_{\partial U_i^\varepsilon} F_i(\phi) \mathbf{n}.
\end{eqnarray*}
The boundary of $U_i^\varepsilon$ consists then of the connected pieces $C_1$ and $C_2$ of level sets of $\phi$,  and two segments, which are $S_1:=U_{i,0}^\varepsilon\cap\{x'_1=0\}$ and $S_2:=U_{i,0}^\varepsilon\cap\{x'_1=|\mathbf{a}|\}$. We have $S_1 = S_2 + \mathbf{a}$. By periodicity of $\phi$, and the fact that the outward unit normal vector on $S_1$ is the opposite of the outward unit normal vector on $S_2$ we get
$$
\int_{S_1} F_i(\phi) \mathbf{n} + \int_{S_2} F_i(\phi) \mathbf{n} = 0.
$$
Moreover, we have
$$
\int_{C_1} F_i(\phi) \mathbf{n} = F_i(\phi(C_1)) \int_{C_1} \mathbf{n} = F_i(\phi(C_1)) R \mathbf{a}.
$$
where $R$ is still the matrix of a rotation of angle $\pi/2$. Besides,
$$
\int_{C_2} F_i(\phi) \mathbf{n} = F_i(\phi(C_2)) \int_{C_2} \mathbf{n} = -F_i(\phi(C_2)) R \mathbf{a}.
$$
It suffices now to consider a function $F_i$ defined by $F_i' = \eta_i^2$ which is not constant on $U_i^\varepsilon$, so we just need to consider a function $\eta_i$ which has compact support in $\phi(U_i^\varepsilon)$, but is not identically zero. This way we get $F_i(\phi(C_2)) \neq F_i(\phi(C_1))$ and we obtain
$$
\int_{U_i^\varepsilon} \nabla F_i(\phi) \neq 0.
$$
We set then
$$
w_0 = \eta_i(\phi) \text{ on } U_i^\varepsilon, \text{ and } 0 \text{ otherwise}.
$$
The function $w_0$ obviously belongs to $\mathcal{I}$, and using (\ref{split}), all the terms in the sum vanish except for the integral on $\hat U_i$, so
$$
\int_C q w_0^2 \neq 0.
$$
This proves (ii).

In the last part of the theorem, we need to prove that whenever $\int_C q w_0^2 \neq 0$, it is proportional to $\mathbf{a}$, where $\mathbf{a}$ is such that all the unbounded periodic trajectories of $q$ in $\Omega$ are $\mathbf{a}-$periodic.

For that purpose, we return  to the previous computations. We know that for $\varepsilon>0$ sufficiently small we have
$$
\int_{\hat U_i^\varepsilon} (\nabla \hat \phi) \hat w^2 = (F_i(\phi(C_1))-F_i(\phi(C_2))) R \mathbf{a}.
$$
Hence, for any $\varepsilon>0$ we have
$$
\mathbf{a} \cdot \int_{\hat U_i^\varepsilon} (\nabla \hat \phi) \hat w^2 = 0,
$$
which remains true at the limit $\varepsilon \to 0$. Using (\ref{split}) we get
$$
\mathbf{a} \cdot \int_C \nabla \phi w^2 = 0,
$$
which gives
$$
R \mathbf{a} \cdot \int_C qw^2 = 0,
$$
for any $w \in \mathcal{I}.$ This is equivalent to say
$$
\int_C qw^2 \in \mathbb{R}\mathbf{a}.
$$
\hfill$\Box$

\end{document}